\RequirePackage{ifthen}
\newboolean{MPA}
\setboolean{MPA}{false}

\ifthenelse {\boolean{MPA}}
{
\documentclass{svjour3}
\smartqed
} {
\documentclass[11pt]{article}
\usepackage[margin=1in]{geometry}
}

\usepackage{amsmath}
\usepackage{graphicx}
\usepackage{color}
\usepackage{amsfonts}
\usepackage{hyperref}
\usepackage{soul}
\usepackage{cancel}
\usepackage{enumitem}

\newcommand{\R}{\mathbb R}
\newcommand{\Z}{\mathbb Z}

\newcommand{\opt}{\textnormal{opt}}
\newcommand{\stt}{\textnormal{s.t.}}

\newcommand{\norminfL}[1]{\left\lVert#1\right\rVert_\infty}
\newcommand{\norminf}[1]{\lVert#1\rVert_\infty}
\newcommand{\absL}[1]{\left\lvert#1\right\rvert}
\newcommand{\abs}[1]{\lvert#1\rvert}
\newcommand{\floor}[1]{\left\lfloor#1\right\rfloor}
\newcommand{\ceil}[1]{\left\lceil#1\right\rceil}
\newcommand{\pare}[1]{\left(#1\right)}
\newcommand{\bra}[1]{\left\{#1\right\}}

\newcommand{\cone}{\operatorname{cone}}

\newcommand{\transp}{\mathsf T}

\usepackage{color}

\newenvironment{prfc}[1][]
{\begin{proof}}
{\end{proof}}

\ifthenelse {\boolean{MPA}}
{

\newenvironment{prf}[1][]
{\begin{proof}}
{\qed \end{proof}}

\newcommand{\qedhere}{\tag*{\qed}}

\newcounter{claim}
\renewenvironment{claim}
{\refstepcounter{claim} \begin{trivlist} \item[] {\bf Claim~\theclaim}\space \itshape}
{\end{trivlist}}

\journalname{Mathematical Programming A}

} {

\usepackage{amsthm}

\newtheorem{theorem}{Theorem}
\newtheorem{lemma}{Lemma}
\newtheorem{claim}{Claim}
\newtheorem{proposition}{Proposition}
\newtheorem{definition}{Definition}
\newtheorem{example}{Example}

\newenvironment{prf}[1][]
{\begin{proof}}
{\end{proof}}

}

\newtheorem{observation}{Observation}

\begin{document}

\title{Proximity in Concave Integer Quadratic Programming \thanks{This work is supported by ONR grant N00014-19-1-2322. Any opinions, findings, and conclusions or recommendations expressed in this material are those of the authors and do not necessarily reflect the views of the Office of Naval Research.}}

\ifthenelse {\boolean{MPA}}
{
\titlerunning{Proximity in Concave Integer Quadratic Programming}
\authorrunning{Alberto Del Pia, Mingchen Ma}

\author{Alberto Del Pia \and Mingchen Ma}
\institute{Alberto~Del~Pia \at
              Department of Industrial and Systems Engineering 
              \& Wisconsin Institute for Discovery \\
              University of Wisconsin-Madison, Madison, WI, USA \\
              \email{delpia@wisc.edu}
              \and
              Mingchen~Ma \at
              Department of Computer Sciences \\
              University of Wisconsin-Madison, Madison, WI, USA \\
              \email{mma54@wisc.edu}}
}
{
\author{Alberto Del Pia
\thanks{Department of Industrial and Systems Engineering \& Wisconsin Institute for Discovery,
             University of Wisconsin-Madison, Madison, WI, USA.
             E-mail: {\tt delpia@wisc.edu}.}
\and
Mingchen Ma
\thanks{Department of Computer Sciences,
             University of Wisconsin-Madison, Madison, WI, USA.
             E-mail: {\tt mma54@wisc.edu}.}
             }
}

\date{\today}

\maketitle

\begin{abstract}
A classic result by Cook, Gerards, Schrijver, and Tardos provides an upper bound of $n \Delta$ on the proximity of optimal solutions of an Integer Linear Programming problem and its standard linear relaxation.
In this bound, $n$ is the number of variables and $\Delta$ denotes the maximum of the absolute values of the subdeterminants of the constraint matrix.
Hochbaum and Shanthikumar, and Werman and Magagnosc showed that the same upper bound is  valid if a more general 
convex function is minimized, instead of a linear function.
No proximity result of this type is known when the objective function is nonconvex.
In fact, 
if we minimize a 
concave quadratic,
no upper bound can be given as a function of $n$ and $\Delta$.
Our key observation is that, in this setting, proximity phenomena still occur, but only if we consider also approximate solutions instead of optimal solutions only.
In our main result
we provide 
upper bounds on the distance between approximate (resp., optimal) solutions to a 
Concave Integer Quadratic Programming problem and optimal (resp., approximate) solutions of its continuous relaxation.
Our bounds are functions of $n, \Delta$, and a parameter $\epsilon$ that controls the quality of the approximation. 
Furthermore, we discuss how far from optimal are our proximity bounds.

\ifthenelse {\boolean{MPA}}
{
\keywords{integer quadratic programming \and quadratic programming \and concave minimization \and proximity \and sensitivity \and subdeterminants}
\subclass{MSC 90C10 \and 90C20 \and 90C26 \and 90C31}
} {}
\end{abstract}

\ifthenelse {\boolean{MPA}}
{}{
\emph{Key words:} integer quadratic programming; quadratic programming; concave minimization; proximity; sensitivity; subdeterminants
}

\section{Introduction}

The relationship between an Integer Linear Programming problem and its standard linear relaxation plays a crucial role in many theoretical and computational aspects of the field, including perfect formulations, cutting planes, and branch-and-bound.
Proximity results study one of the most fundamental questions regarding this relationship: 
Is it possible to bound the distance between optimal solutions to an Integer Linear Programming problem and its standard linear relaxation?
A classic result by Cook, Gerards, Schrijver, and Tardos~\cite{CooGerSchTar86} 
provides the upper bound $n \Delta$ for this distance, where $n$ is the number of variables and $\Delta$ denotes the maximum of the absolute values of the 
subdeterminants of the constraint matrix.
This bound has been recently extended to the mixed-integer case by Paat et al.~\cite{PaaWeiWel18} to $p \Delta$, where $p$ is the number of integer variables.
For other recent proximity results in Integer Linear Programming, we refer the reader to \cite{EisWei18,XuLee19,AliHenOer19}.

Granot and Skorin-Kapov~\cite{GraSko90} show that the upper bound $n \Delta$ is still valid if we minimize a separable convex quadratic objective function over the integer points in a polyhedron.
This result has been further extended to separable convex objective functions by Hochbaum and Shanthikumar~\cite{HocSha90} and by Werman and Magagnosc \cite{WerMag91}.

All the above results feature a convex objective function to be minimized.
Therefore, a natural question is whether proximity phenomena only occur in the presence of convexity.
%
The next example seems to indicate that this is indeed the case.
In fact it shows that, 
with a concave objective,
the distance between optimal solutions of 
the discrete and continuous problems
cannot be bounded by any function of $n, \Delta$.

\begin{example}
\label{ex no bound}
Consider the following 
optimization problem for every $t \in \Z$ with $t \ge 0$:
\begin{align}
\label{pr ex}
\begin{split}
	\min \ & 
	- \pare{x-\frac{1}{4}}^2 \\
	\stt \ & -t\le x \le t+ \frac{3}{4}\\
	&x\in \Z.
\end{split}
\end{align}
Note that all these problems have dimension one ($n=1$) and $\Delta = 1$. 
Clearly, the unique optimal solution to \eqref{pr ex} is $x^d := -t.$ 
If we drop the integer constraint, then the unique optimal solution is 
$x^c := t + \frac{3}{4}.$ 
We have 
$\abs{x^d-x^c} = 2t+\frac{3}{4},$ 
which goes to infinity as $t$ approaches infinity. $\hfill \diamond$
\end{example}

Example~\ref{ex no bound} explains the lack of proximity results in the nonconvex setting.
However, a key observation 
is that the solution $x^* = t$, while not optimal to \eqref{pr ex}, is `almost' optimal.
Furthermore, its distance from $x^c$ is always $\frac 34$.
This simple observation leads us to the question that is at the basis of this work:
Is it possible
to bound the distance between \emph{approximate} (resp., \emph{optimal}) solutions to 
a Nonconvex Integer Programming problem 
and \emph{optimal} (resp., \emph{approximate}) solutions of its continuous relaxation? 
This paper provides the first answers to the posed question.

The optimization problem in Example~\ref{ex no bound} belongs to perhaps the simplest class of Nonconvex Integer Programming problems, formed by Integer Quadratic Programming problems with separable concave objective functions.
Therefore, in this paper we focus on this class of optimization problems.
For these problems, we answer our question in the affirmative and provide explicit upper bounds.
Our bounds are functions of $n, \Delta$, and a parameter $\epsilon$ that controls the quality of the approximation.
Furthermore, we discuss how far from optimal are our proximity bounds.

In the remainder of this section we formally introduce Separable Concave Integer Quadratic Programming and $\epsilon$-approximate solutions. 
With the notation in place, we then formally state our results.

\subsection{Separable Concave Integer Quadratic Programming}

In this paper we denote by \eqref{pr IQP} the Separable Concave Integer Quadratic Programming problem
	\begin{align}
	\label{pr IQP}
	\tag{IQP}
	\begin{split}
	\min \ & \sum_{i=1}^k -q_ix_i^2+h^\transp x \\
	\stt \ & Ax\le b \\
	& x \in \Z^n.
	\end{split}
	\end{align}
	In this formulation we assume that $q_i > 0$ for every $i = 1,\dots, k$.
	Furthermore, we assume that the matrix $A$ is integer, while the remaining data is real.
	Clearly \eqref{pr IQP} subsumes Integer Linear Programming, which can be obtained by setting $k=0$.
We refer the reader to \cite{dPWei14,dP16,dPDeyMol17,dP18,dP20} for recent theoretical results on \eqref{pr IQP}.

	We denote by \eqref{pr QP} the Separable Concave Quadratic Programming problem obtained from \eqref{pr IQP} by dropping the integer constraints, i.e.,
	\begin{align}
	\label{pr QP}
	\tag{QP}
	\begin{split}
	\min \ & \sum_{i=1}^k -q_ix_i^2+h^\transp x \\
	\stt \ & Ax\le b \\
	& x \in \R^n.
	\end{split}
	\end{align}

Throughout this paper, we denote by 
$f(x) := \sum_{i=1}^k -q_ix_i^2+h^\transp x$
the objective function of \eqref{pr IQP} and \eqref{pr QP}, which is quadratic, concave, and separable.
	Furthermore, we let $P$ be the polyhedron defined by $P := \{x \in \R^n \mid Ax\le b\},$ 
	and we denote by $\Delta$ the largest absolute value of the subdeterminants of $A$.

\subsection{$\epsilon$-approximate solution}

In order to state our proximity results, we give the definition of $\epsilon$-approximate solution.
Consider an instance of an optimization problem of the form $\min \{ f(x) \mid x \in S\}$, where $S \subseteq \R^n$.
We assume that this problem has an optimal solution, and we denote it by $x^\opt$.
Let $f_{\max}$ be the maximum value of $f(x)$ on the feasible region $S$.
For $\epsilon \in [0,1]$, we say that a feasible point $x^*$ is an \emph{$\epsilon$-approximate solution} if
\begin{equation*}
f(x^*) - f(x^\opt) \le \epsilon \cdot \pare{f_{\max} - f(x^\opt)}.
\end{equation*}
An intuitive way to interpret this definition is as follows:
If we let $[\alpha,\beta]$ be the smallest interval containing the image of $S$ under $f$, then $f(x^*)$ should lie in the interval $[\alpha, \alpha + \epsilon(\beta - \alpha)]$.
Observe that any feasible point is a $1$-approximation, and only an optimal solution is a $0$-approximation. 
If $f(x)$ has no upper bound on the feasible region, our definition loses its value because any feasible point is an $\epsilon$-approximation for any $\epsilon > 0$. 
Our definition of approximation has been used in earlier works, and we refer to  \cite{NemYud83,Vav92c,BelRog95,KleLauPar06} for more details.

In this work we consider $\epsilon$-approximate solutions to \eqref{pr IQP} and to \eqref{pr QP}.
Clearly, the optimal solution $x^\opt$ and the quantity $f_{\max}$ in the definition of $\epsilon$-approximate solution differ for the two problems because the feasible regions are different.
To avoid confusion, throughout the paper we denote by $x^d$ an optimal solution to \eqref{pr IQP} and by $x^c$ an optimal solution to \eqref{pr QP}.
Similarly, we denote by $f^d_{\max}$ the value $f_{\max}$ in the definition of $\epsilon$-approximate solution to \eqref{pr IQP} and by $f^c_{\max}$ the value $f_{\max}$ in the definition of $\epsilon$-approximate solution to \eqref{pr QP}.

The definition of $\epsilon$-approximate solution is natural for these general problems, and has several useful properties. 
It is well known that, for continuous optimization problems, the definition is insensitive to translations or dilations of the objective function, and that it is preserved under affine linear transformations of the problem.
Similar invariance properties 
hold for discrete optimization problem, 
and are formalized in Lemma~\ref{lem tans} in Section~\ref{sec lemmas}.

\subsection{Our results}


We are ready to state our proximity result for Separable Concave Integer Quadratic Programming.

\begin{theorem}
\label{th main}
Consider a problem \eqref{pr IQP}, and the corresponding continuous problem \eqref{pr QP}. 
Suppose that both problems 
have an optimal solution.
Then:
\begin{enumerate}[label={(\roman*)}] 
\item
\label{th main 1}
Let $x^c$ be an optimal solution to \eqref{pr QP}.
Then, $\forall \epsilon \in (0,1]$, there is 
an $\epsilon$-approximate solution $x^*$ to \eqref{pr IQP}
such that 
$\norminf{x^c-x^*} \le n\Delta \pare{\frac{10\Delta}{\epsilon}+1}^k.$
\item
\label{th main 2}
Let $x^d$ be an optimal solution to \eqref{pr IQP}.
Then, $\forall \epsilon \in (0,1]$, there is an $\epsilon$-approximate solution $x^\star$ to \eqref{pr QP} such that 
$\norminf{x^\star-x^d} \le n\Delta \pare{\frac{10\Delta}{\epsilon}+1}^k.$
\end{enumerate}
\end{theorem}

In particular, note that the bounds in Theorem~\ref{th main} do not depend on the right-hand side vector $b$ in \eqref{pr IQP} and \eqref{pr QP}.
The proof of Theorem~\ref{th main} is given in Section~\ref{sec proof}.
Since $k \le n$, Theorem~\ref{th main} implies that, for every optimal solution to one of the two problems, and for every $\epsilon \in (0,1]$, there is an $\epsilon$-approximate solution to the other problem at distance bounded 
by a function of $n,\Delta,\epsilon$.
In particular, this distance is independent on the objective function and on the vector $b$.
Note that, for $k=0$, problem \eqref{pr IQP} is an Integer Linear Programming problem, while \eqref{pr QP} is its standard linear relaxation.
In this setting, our bounds in Theorem~\ref{th main} reduce to $n\Delta$ for every $\epsilon > 0$.
Therefore, the proximity bound by Cook et al.~\cite{CooGerSchTar86} can be obtained as a corollary to Theorem~\ref{th main}.

In Section~\ref{sec tight}, we discuss how far from optimal are our upper bounds in Theorem~\ref{th main}.
At the heart of our tightness results lies a special polytope, denoted by $\bar P$, and which is used with several different objective functions.
In particular, using the notation of Theorem~\ref{th main}, we show that any upper bound on $\norminf{x^c-x^*}$ or on $\norminf{x^\star-x^d}$ must grow at least linearly with $\frac 1 \epsilon,$ $n$, and $\Delta$.
Furthermore, we show that the neighborhood of $x^c$ considered by Cook et al.~\cite{CooGerSchTar86}, namely $\{x \in P \cap \Z^n \mid \norminf{x^c-x} \le n\Delta\}$ might contain only arbitrarily bad solutions to \eqref{pr IQP}, i.e., vectors that are not $\epsilon$-approximate solution to \eqref{pr IQP}, for any $\epsilon \in (0,1)$.
The polytope $\bar P$ also allows us to show that the Integer Linear Programming bound $n \Delta$ by Cook et al. 
is best possible.
To the best of our knowledge this tightness result was known only for $\Delta = 1$ 
(see 
 page 241 in \cite{SchBookIP}).

\section{Three simple lemmas}
\label{sec lemmas}

In this section we present three lemmas that will be used in the proof of Theorem~\ref{th main}.

Our first lemma formalizes the invariance properties of $\epsilon$-approximate solutions to optimization problems with integer constraints.
The proof is standard.
This result will allow us to greatly simplify the notation in the main proof. 

\begin{lemma}
\label{lem tans}
    Consider an optimization problem of the form
	\begin{align}
	\label{pr O}
	\tag{O}
	\begin{split}
	\min \ & f(x) \\
	\stt \ & x \in S \cap \Z^n,
	\end{split}
	\end{align}
	where $S \subseteq \R^n$.
	Let $M \in \Z^{n\times n}$ be a unimodular matrix, and let $t\in \Z^n$.
	For any $\alpha, \beta \in \R$ with $\alpha > 0$, consider the optimization problem
	\begin{align}
	\label{pr O'}
	\tag{O'}
	\begin{split}
	\min \ & \alpha f(M^{-1} (y-t)) + \beta \\
	\stt \ & y\in U(S) \cap \Z^n,
	\end{split}
	\end{align}
	where $U(S) := \{y \in \R^n \mid y=Mx+t, \ x\in S\}$.
	Then, for every $\epsilon$-approximate solution to $\eqref{pr O}$, denoted by $x^*$, the vector $Mx^* + t$ is an $\epsilon$-approximate solution to $\eqref{pr O'}$.
	Viceversa, for every $\epsilon$-approximate solution to $\eqref{pr O'}$, denoted by $y^*$, the vector $M^{-1}(y^*-t)$ is an $\epsilon$-approximate solution to $\eqref{pr O}$.
\end{lemma}	

\begin{prf}
We prove the first statement of the lemma, the second one being symmetric.
Let $x^*$ be an $\epsilon$-approximate solution to $\eqref{pr O}$.
We show that the vector $y^* := Mx^* + t$ is an $\epsilon$-approximate solution to $\eqref{pr O'}$.

Since $M$ and $t$ are integer, for every feasible solution $x$ to $\eqref{pr O}$, the vector $y = Mx+t$ is feasible to $\eqref{pr O'}$.
Viceversa, since $M^{-1}$ and $t$ are integer, for every feasible solution $y$ to $\eqref{pr O'}$, the vector $x = M^{-1}(y-t)$ is feasible to $\eqref{pr O}$.
In both cases, the relation between the cost of $x$ and $y$ is given by $g(y) = \alpha f(x) + \beta$,
where $g(y) := \alpha f(M^{-1} (y-t)) + \beta$ denotes the objective function of \eqref{pr O'}.

Let $x^d$ be an optimal solution to $\eqref{pr O}$, and let $y^d$ be an optimal solution to $\eqref{pr O'}$.
Furthermore, let $f_{\max}$ be the maximum value of $f(x)$ on the feasible region of $\eqref{pr O}$, and let $g_{\max}$ be the maximum value of $g(y)$ on the feasible region of $\eqref{pr O'}$.
Since $\alpha > 0$, the above argument in particular implies $g(y^d) = \alpha f(x^d) + \beta$, and $g_{\max} = \alpha f_{\max} + \beta$.
If $g_{\max} = g(y^d)$, then $y^*$ is an optimal solution to $\eqref{pr O'}$ and we are done. Otherwise, we have
\begin{equation*}
    \frac{g(y^*)-g(y^d)}{g_{\max}-g(y^d)}
    =\frac{\pare{\alpha f(x^*)+\beta}-\pare{\alpha f(x^d)+\beta}}{\pare{\alpha f_{\max}+\beta}-\pare{\alpha f(x^d)+\beta}}
    =\frac{f(x^*)-f(x^d)}{f_{\max}-f(x^d)}\le \epsilon.
\end{equation*}
Thus $y^*$ is an $\epsilon$-approximate solution to \eqref{pr O'}.
\end{prf}

Next, we define a polyhedral cone which will be heavily used in the proof of Theorem~\ref{th main}, and we present some of its properties.
We remark that this cone has been used in several papers to obtain proximity results, including \cite{CooGerSchTar86,GraSko90,HocSha90,dP20}.
Let $A$ be a matrix with $n$ columns and let $x^a,x^b \in \R^n$.
Let $A_1$ be the matrix that contains all rows $u$ in $A$ for which $ux^a \le ux^b$. 
Similarly, let $ A_2$ be the matrix that contains all rows $u$ in $A$ for which $ ux^a \ge ux^b$. 
We define the polyhedral cone 
\begin{align*}
T(A, x^a, x^b) := \bra{x \in \R^n \mid A_1 x \le 0, \ A_2 x \ge 0}.
\end{align*}
From the definition of the cone, we obtain $x^a-x^b \in T(A, x^a, x^b).$

The next lemma is well-known, see, e.g., \cite{CooGerSchTar86}.
Since we are unable to find a complete proof in the literature, we present it here.
\begin{lemma} 
\label{lem 1}
Let $A$ be an integer matrix with $n$ columns, let $\Delta$ be the largest absolute value of the 
subdeterminants of $A$, and let $x^a,x^b \in \R^n$.
Then there exists a finite subset $V$ of $\Z^n$ such that $T(A, x^a, x^b)  = \cone V,$ and for every $v \in V$, we have $\norminf{v} \le \Delta$.
\end{lemma}

\begin{prf}
Let $T := T(A, x^a, x^b).$ 
We partition $T$ into pointed polyhedral cones by intersecting it with the $2^n$ orthants of $\R^n$, which we denote by $O_1, \dots, O_{2^n}$.
Namely, we let $T_i := T \cap O_i$, for $i=1,\dots,2^n$, and observe that $T = \bigcup_{i=1}^{2^n} T_i.$ 
In order to prove the lemma, it suffices to show that, for every $i=1,\dots,2^n$, there exists a finite subset $V_i$ of $\Z^n$ such that $T_i  = \cone V_i,$ and for every $v \in V_i$, we have $\norminf{v} \le \Delta$.
This is because the set $V = \bigcup_{i=1}^{2^n} V_i$ then satisfies the thesis of the lemma.

Let us now consider a single $T_i$, for some $i \in \{1,\dots,2^n\}$.
We assume that $T_i$ arises from the intersection of $T$ with the nonnegative orthant, i.e., $T_i = \{x \in \R^n \mid x \in T, \ x \ge 0\}$, the other cases being symmetric.
The set $T_i$ is a pointed polyhedral cone.
Since $A$ is integer, $T_i$ is a rational cone.
Therefore, there exists a finite set of vectors $V_i = \{r^1,\dots,r^m\} \subset \R^n$ such that $T_i = \cone V_i.$ 
Here we can assume that for every $j=1,\dots,m,$ the vector $r^j$ is not a proper conic combination of other vectors in $V_i$, that is to say, each $r^j$ is an extreme ray of $T_i.$ 

Let us now consider a single vector $r^j$, for some $j \in \{1,\dots,m\}$.
We now show that we can scale $r^j$ so that it is integer and with infinity norm at most $\Delta.$
From Theorem 3.35 in \cite{ConCorZamBook}, we know that $r^j$ satisfies at equality $n-1$ linearly independent inequalities in the system $A_1 x \le 0$, $A_2 x \ge 0$, $x \ge 0$.
Let $e_k$ be a vector of the standard basis of $\R^n$ that is linearly independent from all the rows of $A_1$, $A_2$, and the identity matrix which correspond to the $n-1$ linearly independent inequalities.
Note that we have $r^j_k \neq 0,$ since otherwise we obtain $r^j = 0$, which contradicts the fact that $r^j$ is an extreme ray of $T_i.$ 
Since $T_i$ is contained in the nonnegative orthant, we have $r^j_k>0.$
Denote by $Dx=e_1$ the system of equations containing $x_k=1$ and the $n-1$ equations arising by setting to equality the $n-1$ linearly independent inequalities discussed above, where $e_1$ is the first vector of the standard basis of $\R^n$.
Note that the matrix $D$ is invertible.
The vector $r := D^{-1} e_1$ is a solution to the system, and is a scaling of the vector $r^j$.
Note that each entry of $r$ coincides with an entry of the matrix $D^{-1}$.
By Cramer's rule, each entry of $D^{-1}$ is a fraction with denominator $\det(D)$ and numerator with absolute value at most $\Delta.$
Thus, the vector $\abs{\det(D)} \cdot r$ is a scaling of $r^j$ that is integer and with $\norminf{\abs{\det(D)} \cdot r} \le \Delta.$ 
Hence, we can assume that each vector $r^j$ is integer and with infinity norm at most $\Delta.$
\end{prf}


The next lemma will often be used in the proof of Theorem~\ref{th main} to show that a given vector is in our polyhedron $P$.

\begin{lemma} \label{lem 2}
Let $P = \{x\in \R^n \mid A x \le b\}$ be a polyhedron,
let $x^a,x^b \in P$.
Let $x^\circ$ be a vector in $\R^n$ that can be written in the following two ways:
    \begin{align*}
        x^\circ=x^1+\sum_{i=1}^m \alpha_i v^i, \qquad 
        x^\circ=x^2-\sum_{i=1}^m \beta_i v^i,
    \end{align*}
    where $x^1, x^2 \in P$
    and, for $i =1,\dots,m$, $\alpha_i, \beta_i,$ are nonnegative numbers and $v^i \in T(A, x^a, x^b)$.
    Then $x^\circ \in P.$ 
\end{lemma}

\begin{prf}
Let $A_1, A_2$ from the definition of $T(A, x^a, x^b)$ and $b_1, b_2$ the corresponding sub-vectors of $b.$
Since $x^1,x^2 \in P,$ we obtain
\begin{align*}
    & A_1 x^\circ
    = A_1 \pare{x^1+\sum_{i=1}^m \alpha_i v^i}
    \le b_1+\sum_{i=1}^m \alpha_i A_1 v^i
    \le b_1\\
    & A_2 x^\circ
    = A_2 \pare{x^2-\sum_{i=1}^m \beta_i v^i}
    \le b_2-\sum_{i=1}^m \beta_i A_2 v^i
    \le b_2,
\end{align*}
where the last inequalities follow because $A_1 v^i \le 0$ and $A_2 v^i \ge 0$ from the definition of $T(A, x^a, x^b)$.
This implies that $Ax^\circ \le b,$ hence $x^\circ \in P.$ 
\end{prf}


We are now ready to present our proof of Theorem~\ref{th main}.

\section{Proof of Theorem~\ref{th main}}
\label{sec proof}

Let $x^c$ be an optimal solution to \eqref{pr QP}, let $x^d$ be an optimal solution to \eqref{pr IQP}, and let $\epsilon >0$. In this section, we present our proof of Theorem~\ref{th main}. To do so, we will construct an $\epsilon$-approximate solution $ x^* $ to \eqref{pr IQP} and an $\epsilon$-approximate solution $ x^\star $ to \eqref{pr QP}.

We first give a brief outline of our proof. 
	In Section~\ref{sec algorithm}, we design a recursive algorithm which takes in input $x^c, x^d, P$ and that outputs a point $ x^\ell \in P $. In the algorithm, we use several times cones of the form $T(A, \cdot ,\cdot)$ to construct a path inside $ P $, which starts at $ x^c $, ends at $ x^\ell $ and contains at most $ k+1 $ points. The special structure of this path enables us to upper bound $ \norminf{x^c-x^\ell} $ by a function of $ n,\Delta,k,\epsilon.$ 
	In Section~\ref{sec property}, we study some properties of $x^\ell$ and we consider separately two cases. In the first case, $ \norminf{x^\ell-x^d} $ can be bounded by a function of $ n,\Delta,k,\epsilon.$ In this case, we can then also bound $ \norminf{x^c-x^d} $ by a function of $ n,\Delta,k,\epsilon.$ As a consequence, we can conclude the proof in the first case by choosing $ x^* $ to be $ x^d $ and by choosing $ x^\star $ to be $ x^c $. 
	In the second case, $ \abs{x^\ell_i-x^d_i} $ is large for every index $ i \in \{1,\dots,k \}$ such that $ x^\ell_i-x^d_i \neq 0 .$ In this case, in Section~\ref{sec integral point} we use vectors $x^\ell,x^d$ to construct an integer vector $x^*$ which is close to $x^\ell$, that is, $\norminf{x^\ell-x^*} \le n\Delta.$ This in particular implies that $\norminf{x^c-x^*}$ can be bounded by a function of $ n,\Delta,k,\epsilon.$ In Section~\ref{sec property 2} we further study the vector $x^*$. The properties obtained allow us to prove, in Section~\ref{sec analysis 1}, that $x^*$ is an $\epsilon$-approximate solution to \eqref{pr IQP}. This is done by providing an upper bound on $ f(x^*)-f(x^d) $ and a lower bound on $ f^d_{\max}-f(x^d) .$ In Section~\ref{sec star} we define a vector $ x^\star $, based on $ x^* $, with the property that $ \norminf{x^\star-x^d}=\norminf{x^c-x^*} $. Next, in Section~\ref{sec analysis 2}, we prove that $ x^\star $ is an $ \epsilon$-approximate solution to \eqref{pr QP}. This concludes the proof in the second case, and our outline of the proof of Theorem~\ref{th main}. We are now ready to present the full proof.

In order to simplify the notation in the remainder of the proof, in the next claim we employ Lemma~\ref{lem tans}.

\begin{claim}
\label{claim eqt}
    We can assume without loss of generality that $x^d$ is the origin and that $f(x^d)=0$.
\end{claim}

\begin{prf}
We apply Lemma~\ref{lem tans} as follows:
Problem \eqref{pr O} is \eqref{pr IQP}, the matrix $M$ is the identity, $t := -x_d$, $\alpha := 1$, and $\beta := -f(x^d)$.
Problem \eqref{pr O'} in Lemma~\ref{lem tans} then takes the form
	\begin{align}
	\label{pr O' inproof}
	\begin{split}
	\min \ & f(y+x^d)-f(x^d) \\
	\stt \ & Ay\le b-Ax^d\\
	& y \in \Z^n.
	\end{split}
	\end{align}
	The objective function of \eqref{pr O' inproof} can then be explicitly written as 
	$\sum_{i=1}^k - q_i y_i^2 + {h'}^\transp y,$
	where the vector $h'$ is defined by $h'_i := h_i - 2q_i x_i^d$ for $i=1,\dots,k$ and $h'_i := h_i$ for $i=k+1,\dots,n$.
	In particular, the coefficients $q_i$ of the quadratic monomials are identical in \eqref{pr IQP} and in \eqref{pr O' inproof}. 
	Furthermore, note that the constraint matrix of \eqref{pr O' inproof} is the same of \eqref{pr IQP} and so the two problems have the same $\Delta$.

The transformation that maps problem \eqref{pr IQP} into \eqref{pr O' inproof} is $y = x - x_d$.
Consider vectors $\bar x, \tilde x, \bar y, \tilde y \in \R^n$ such that $\bar y = \bar x - x_d$ and $\tilde y = \tilde x - x_d$.
Then distances are maintained since $\bar y - \tilde y = \bar x - \tilde x$.

Lemma~\ref{lem tans} implies that $\epsilon$-approximate solutions are mapped to $\epsilon$-approximate solutions, thus, without loss of generality we can consider problem \eqref{pr O' inproof} instead of \eqref{pr IQP}.
The optimal solution $x^d$ of \eqref{pr IQP} is mapped to the origin, which is then an optimal solution of \eqref{pr O' inproof}, and the optimal cost is zero.
\end{prf}


\subsection{Construction of the vector $x^\ell$}\label{sec algorithm}

This section of the proof is devoted to the construction of a special vector in $P$ that we denote by $x^\ell$. 
The vector $x^\ell$ is obtained via a recursive algorithm which utilizes the vectors $x^c$ and $x^d.$
To begin with, we introduce a claim which will play a key role in the iterative step of our algorithm.

\begin{claim}
\label{claim onestep}
Let $x^a \in P$ and let 
$Z := \{i \in \{1,\dots, k\} \mid x^a_i=0\}$. 
Assume $Z \neq \{1,\dots,k\}$,
and let 
$s$ be an index such that $\abs{x^a_s}=\min \{\abs{x^a_i} \mid i \in \{1,\dots,k\} \setminus Z\}.$
Assume that 
$\norminf{x^a} > \Delta \abs{x^a_s}$.
Then there exists a vector $x^b \in P$ such that
\begin{enumerate}[label={(\roman*)}, leftmargin=*] 
\item
\label{cond onestep1}
$x^b_i = 0$ for every $i \in Z \cup \{s\}$;
\item
\label{cond onestep2}
$\norminf{x^a-x^b} \le \Delta \abs{x^a_s}$;
\item
\label{cond onestep3}
For $i =1,\dots,m$, there exist nonnegative scalars $\alpha_i,$ $\beta_i$ and vectors $v^i \in T(A,x^a,x^d)$ such that $x^b$ can be expressed in the following two ways:
\begin{align*}
    x^b =x^d+\sum_{i=1}^m \alpha_i v^i, \qquad
    x^b =x^a-\sum_{i=1}^m \beta_i v^i.
\end{align*}
\end{enumerate}
\end{claim}

\begin{prf}
Our fist task is that of defining the vector $x^b$.
Denote by $\tilde A x\le \tilde b$ the system obtained from $Ax\le b$ by adding the inequalities 
$x_i \le 0$, $-x_i \le 0$, for all $i \in Z$.
We remark that the largest absolute value of a 
subdeterminant
of $\tilde A$ is $\Delta$. 
Let $\tilde P := \{x \in \R^n \mid \tilde A x \le \tilde b \}$, and note that the vectors $x^a$ and $x^d$ are in $\tilde P$.
Let $\tilde T := T(\tilde A, x^a, x^d)$.
From Lemma~\ref{lem 1}, applied to $\tilde A$, $x^a$, and $x^d$, we know that there exists a finite subset $\tilde V$ of $\Z^n$ such that $\tilde T = \cone \tilde V,$ and for every $v \in \tilde V$, we have $\norminf{v} \le \Delta$.
Since $x^a - x^d \in \tilde T$, Caratheodory's theorem implies that there exist $m\le n$ vectors $v^1,\dots,v^m \in \tilde V$ and $m$ positive scalars $\alpha_1,\dots, \alpha_m$ such that
\begin{align}
\label{eq op}
x^a = x^a-x^d=\sum_{i=1}^{m}\alpha_iv^i.    
\end{align}

We pick those vectors $v^{i}$ such that $v^{i}_s$ has the same sign as $x^a_s$. Without loss of generality, we can assume that these vectors are $v^{1},\dots,v^{r}$, where $r \le m$.
We now show that there exist nonnegative scalars $\lambda_{1}, \dots, \lambda_{r}$ that satisfy $\lambda_{1}\le \alpha_{1},\dots, \lambda_{r}\le \alpha_{r},$ and such that 
\begin{align}
    \label{eq also this}
x^a_s = \sum_{i=1}^{r}\lambda_{i}v^i_s.
\end{align}
From \eqref{eq op}, we obtain
$x^a_s = \sum_{i=1}^m \alpha_i v^i_s.$
Since for $i=1,\dots,r,$ $v^i_s$ has the same sign as $x^a_s$ and for $i=r+1,\dots,m,$ $v^i_s$ has the opposite sign of $x^a_s$ or $v^i_s=0,$ we have
\begin{align*}
0<\abs{x^a_s}=\abs{\sum_{i=1}^{r}\alpha_{i}v^i_s}-\abs{\sum_{i=r+1}^{m}\alpha_{i}v^i_s} \le \abs{\sum_{i=1}^{r}\alpha_{i}v^i_s}.
\end{align*}
Using continuity, we know that there exist nonnegative scalars $\lambda_{1}\le \alpha_{1},\dots, \lambda_{r}\le \alpha_{r},$ such that
$\abs{x^a_s}=\abs{\sum_{i=1}^{r}\lambda_{i}v^i_s}.$
Since each $v^i_s$ above has the same sign as $x^a_s$ and each $\lambda_i \ge 0,$ we know that \eqref{eq also this} holds.

We are finally ready to define the vector $x^b$ as
\begin{align}
\label{np}
    x^b:=x^a-\sum_{i=1}^{r}\lambda_{i}v^{i}.
\end{align}
From~\eqref{eq op}, we can write $x^b$ in the form
\begin{align}
\label{newnp}
    x^b=
    x^d + \sum_{i=1}^{r}(\alpha_{i}-\lambda_{i})v^{i}+\sum_{i=r+1}^{m}\alpha_{i}v^{i}.
\end{align}
Since $x^a,x^d \in \tilde P$, from Lemma~\ref{lem 2} we know that $x^b$ is in $\tilde P$ as well.
Since $\tilde P \subseteq P$, we obtain that $x^b \in P.$
Next we show that \ref{cond onestep1}, \ref{cond onestep2}, \ref{cond onestep3} hold.

\smallskip

\ref{cond onestep1}. 
Note that $\tilde P$ satisfies equations $x_i=0$, for $i\in Z,$ hence $x^b_i=0$ for every $i \in Z.$
Furthermore, from the definition of $x^b$, and using \eqref{eq also this}, we have
\begin{align*}
x^b_s=x^a_s-\sum_{i=1}^{r}\lambda_{i}v^i_s
=x^a_s-x^a_s=0.
\end{align*}

\smallskip

\ref{cond onestep2}.
From the definition of $x^b$ we have
$\norminf{x^a-x^b}
=\norminf{\sum_{i=1}^{r}\lambda_{i}v^{i}}.$
Denote by $l$ the index such that $\norminf{\sum_{i=1}^{r}\lambda_{i}v^{i}} = \abs{(\sum_{i=1}^{r}\lambda_{i}v^{i})_l}$. 
Then we have 
\begin{align*}
\norminfL{\sum_{i=1}^{r}\lambda_{i}v^{i}}
= \absL{\sum_{i=1}^{r}\lambda_{i}v^{i}_l}
= \absL{\sum_{i=1}^{r}\lambda_{i} \frac{v^{i}_l}{v^{i}_s}v^{i}_s}
\le \sum_{i=1}^{r}\lambda_{i} \absL{\frac{v^{i}_l}{v^{i}_s}} \abs{v^{i}_s}.
\end{align*}    
Since, for every $i=1,\dots,r$, the vector $v^{i}$ is integer and $\norminf{v^{i}} \le \Delta$, we know that $\absL{\frac{v^{i}_l}{v^{i}_s}} \le \Delta.$ 
     Thus, 
     \begin{align*}
     \sum_{i=1}^{r}\lambda_{i} \absL{\frac{v^{i}_l}{v^{i}_s}} \abs{v^i_s} 
     \le \Delta \sum_{i=1}^{r}\lambda_{i}\abs{v^{i}_s}
     = \Delta \absL{\sum_{i=1}^{r}\lambda_{i}v^{i}_s} 
     = \Delta \abs{x^a_s}.
     \end{align*}
     The first equality holds because all $v^{1}_s,\dots,v^{r}_s$ have the same sign, and the last equality follows from \eqref{eq also this}.
     This completes the proof of \ref{cond onestep2}.

\ref{cond onestep3}. 
We notice that $\tilde T \subseteq T(A, x^a, x^d),$ which implies that $v^1, \dots, v^m \in T(A, x^a, x^d).$ 
So \ref{cond onestep3} follows directly from \eqref{np} and \eqref{newnp}.
\end{prf}

We are now ready to state our algorithm that constructs the vector $x^\ell.$
We recursively define a sequence of vectors in $P$ denoted by $x^{0},x^{1},x^{2},\dots$.
The last vector in this sequence is indeed the vector $x^\ell$ that we wish to obtain.
To define this sequence of vectors,
we first recursively define the following $k$ scalars:
\begin{align*}
    \chi_1 &:= \frac{8n\Delta}{\epsilon}+2n\Delta \\
	\chi_j &:= 2n\Delta +\frac{8}{\epsilon} \pare{\sum_{i=1}^{j-1}\Delta\chi_i+n\Delta} && j=2,\dots,k.
\end{align*}

For every vector $x^j$ in the sequence, it will be useful to partition the set $\{1,\dots,k\}$ into the two sets
\begin{align*}
Z^j := \bra{i \in \{1,\dots, k\} \mid x^{j}_i=0}, \qquad 
N^j := \bra{i \in \{1,\dots, k\} \mid x^{j}_i \neq 0}.
\end{align*}
We start the sequence by setting $x^{0} := x^c.$
Now assume that we have constructed the vectors $x^{0},x^1,\dots, x^{j}.$ 
We state the next iteration of the algorithm.
In this iteration, either the algorithm sets $\ell := j$ and terminates, or it  constructs the next vector $x^{j+1}.$ 


If $x^{j}$ satisfies $\abs{x^{j}_i}> \chi_{j+1}$ for every $i \in N^j$ then we set $\ell := j$ and terminate.
Otherwise, we have $N^j \neq \emptyset$ and $\abs{x^{j}_s} \le \chi_{j+1}$, where $s$ is an index such that $\abs{x^{j}_s}=\min \{\abs{x^{j}_i} \mid i \in N^j\}.$

If $\norminf{x^{j}}\le \Delta\abs{x^{j}_s}$, we set $\ell := j$ and terminate.
Otherwise, we have $\norminf{x^{j}}> \Delta\abs{x^{j}_s}.$ 
Then $x^{j}$ satisfies the assumptions of Claim~\ref{claim onestep}.
The next vector $x^{j+1}$ in the sequence is then defined to be the vector $x^b \in P$ in the statement of Claim~\ref{claim onestep}, invoked with $x^a = x^j$.

\smallskip

This concludes our definition of the sequence $x^{0},x^1,\dots, x^{\ell}$ of vectors in $P$.
Note that the sequence contains at most $k+1$ points, i.e., $\ell \le k$.
In fact, according to Claim~\ref{claim onestep}\ref{cond onestep1}, we know that $Z^{j+1}$ has at least one more element than $Z^{j}$, for every $j=0,1,\dots$.

\subsection{Properties of the vector $x^\ell$}\label{sec property}


For ease of notation, we define 
\begin{align*}
\psi_j &:= \sum_{i=1}^j \Delta \chi_i &  j=1,\dots,k,
\end{align*}
and obtain an upper bound on $\psi_k$.

\begin{claim}\label{bound}
We have $\psi_k \le n\Delta(\frac{10\Delta}{\epsilon}+1)^k - n\Delta.$
\end{claim}

\begin{prfc}
The number $\psi_j+n\Delta$ can be upper bounded as follows:
\begin{align*}
\psi_j+n\Delta 
&= 2n\Delta^2+\pare{\frac{8\Delta}{\epsilon}+1}(\psi_{j-1}+n\Delta) \\
& \le\Delta(\psi_{j-1}+n\Delta)+\pare{\frac{8\Delta}{\epsilon}+1}(\psi_{j-1}+n\Delta) \\
& \le \pare{\frac{9\Delta}{\epsilon}+1}(\psi_{j-1}+n\Delta).
\end{align*}
The equality holds by definition of $\psi_j$ and $\chi_j$. 
The first inequality follows from the fact that $n\Delta \le \psi_i$ for every $i=1,\dots,k,$ while the second inequality is correct because $\Delta \le \frac{\Delta}{\epsilon}.$

Then we have
\begin{align*}
\psi_k+n\Delta \le \pare{\frac{9\Delta}{\epsilon}+1}^{k-1}(\psi_1+n\Delta).
\end{align*}
Since
\begin{align*}
\psi_1+n\Delta
=\frac{8n\Delta^2}{\epsilon}+2n\Delta^2+n\Delta
=n\Delta\pare{\frac{8\Delta}{\epsilon}+2\Delta+1} 
\le n\Delta\pare{\frac{10\Delta}{\epsilon}+1},
\end{align*}
we get
\begin{equation*}
\psi_k+n\Delta 
\le n\Delta\pare{\frac{9\Delta}{\epsilon}+1}^{k-1}\pare{\frac{10\Delta}{\epsilon}+1}
\le n\Delta\pare{\frac{10\Delta}{\epsilon}+1}^k. \qedhere
\end{equation*}
\end{prfc}

In the next claim we highlight some properties of $x^\ell$ that will be used later.

\begin{claim}
\label{claim xell}
The vector $x^\ell$ satisfies the following properties:
\begin{enumerate}[label={(\alph*)}] 
\item
\label{xell a}
$x^c-x^\ell \in P;$
\item
\label{xell b}
$\norminf{x^c-x^\ell} \le \psi_\ell;$
\item
\label{xell c}
At least one of the following holds:
\begin{enumerate}[label={(c-\arabic*)}]
\item 
\label{xell c-1}
$\norminf{x^\ell} \le \Delta\chi_{\ell+1},$ with $\ell \le k-1;$
\item 
\label{xell c-2}
$\abs{x^\ell_i} > \chi_{\ell+1}$ for every $i \in N^\ell$.
\end{enumerate}
\end{enumerate}
\end{claim}

\begin{prf}
\ref{xell a}.
We prove the stronger statement that $x^c-x^{j} \in P$ for every $j = 0,\dots,\ell,$ by induction on $j$.
The base case is $j=0$, and it holds since $x^c-x^{0} = 0 = x^d \in P$.

Next, we show the inductive step. 
We assume that the result is true for $j=t,$ and we prove it for $j=t+1.$ 
From our definition of the sequence $x^{0},x^1,\dots, x^{\ell}$, the vector $x^{t+1}$ is obtained from $x^t$ as described in Claim~\ref{claim onestep}, where $x^a = x^{t}$ and $x^b = x^{t+1}$. 
Claim~\ref{claim onestep}\ref{cond onestep3} implies that for $i =1,\dots,m$ there exist nonnegative scalars $\alpha_i,$ $\beta_i$ and vectors $v^i \in T(A, x^t,x^d)$ such that 
\begin{align*}
x^c-x^{t+1} & = (x^c-x^{t})+(x^{t}-x^{t+1})=(x^c-x^{t})+\sum_{i=1}^m \beta_i v^i \\
x^c-x^{t+1} & = x^c-x^d-\sum_{i=1}^m \alpha_i v^i= x^c-\sum_{i=1}^m \alpha_i v^i.
\end{align*}
Clearly $x^c \in P$ and, from the induction hypothesis, $x^c-x^t \in P$ as well.
Then Lemma~\ref{lem 2} implies that $x^c-x^{t+1}\in P.$ 
This concludes our proof that $x^c-x^{j} \in P$ for every $j = 0,\dots,\ell.$
Therefore $x^c-x^\ell \in P$, concluding the proof of \ref{xell a}.

\smallskip

\ref{xell b}. 
According to Claim~\ref{claim onestep}\ref{cond onestep2},
and from the definition of the sequence, 
we know that
\begin{align*}
\norminf{x^{j-1}-x^{j}}\le \Delta\abs{x_s^{j-1}}\le \Delta\chi_j \qquad j=1, \dots, \ell.
\end{align*}
Thus, we have
\begin{align*}
\norminf{x^c-x^\ell} 
\le \sum_{j=1}^{\ell}\norminf{x^{j-1}-x^{j}} 
\le \sum_{j=1}^{\ell}\Delta\chi_j
=\psi_\ell.
\end{align*}

\smallskip

\ref{xell c}.
This proof follows from the definition of $x^\ell.$ 
In fact, since $x^\ell$ is the last point in the sequence, it must satisfy at least one of the two termination conditions.
If $\abs{x^\ell_i} > \chi_{\ell+1}$ for every $i \in N^\ell$, then \ref{xell c-2} holds and we are done.
Note that, if $\ell = k$, then $Z^\ell=\{1,\dots,k\},$ and this termination condition is triggered.
Otherwise, we have $\norminf{x^\ell}\le \Delta\abs{x^\ell_s}$ and $\ell \le k-1.$ 
Observing that $\abs{x^\ell_s} \le \chi_{\ell+1}$ from the construction of the sequence, we obtain \ref{xell c-1}.
\end{prf}

From \ref{xell c}, the vector $x^\ell$ satisfies at least one of the two properties \ref{xell c-1}, \ref{xell c-2}. 
Next we show that, if $x^\ell$ satisfies \ref{xell c-1},
then Theorem~\ref{th main} holds with $x^* = x^d$ and $x^\star = x^c$.
So assume that the vector $x^\ell$ satisfies property \ref{xell c-1}. 
We obtain
\begin{align*}
\norminf{x^c-x^d} 
\le \norminf{x^c-x^\ell}+ \norminf{x^\ell-x^d} 
\le \psi_\ell + \Delta\chi_{\ell+1}
=\psi_{\ell+1}\le \psi_k,
\end{align*}
where the second inequality follows from \ref{xell b} and \ref{xell c-1}. 
Hence the distance between $x^d$ and $x^c$ is upper bounded by $\psi_k$, which is at most $n\Delta(\frac{10\Delta}{\epsilon}+1)^k$ from Claim~\ref{bound}.
As a consequence, in this case, we conclude the proof of Theorem~\ref{th main}\ref{th main 1} with $x^* = x^d$ and of Theorem~\ref{th main}\ref{th main 2} with $x^\star = x^c$.
Therefore, in the remainder of the proof, we assume that $x^\ell$ satisfies \ref{xell c-2}.

\subsection{Construction of the vector $x^*$}\label{sec integral point}

This section of the proof is devoted to the construction of the vector $x^*$ in the statement of Theorem~\ref{th main}\ref{th main 1}.
In particular, $x^*$ lies in a neighborhood of the vector $x^\ell$.

Denote by $\bar A x\le \bar b $ the system obtained from $Ax\le b$ by adding the inequalities 
$x_i \le 0$, $-x_i \le 0$, for all $i \in Z^\ell$.
Note that the largest absolute value of a subdeterminant of $\bar A$ is $\Delta$. 
Let 
\begin{align*}
\bar P := \bra{x \in \R^n \mid \bar A x \le \bar b }.
\end{align*}
Note that $\bar P \subseteq P$ and that the vectors $x^\ell$ and $x^d$ are in $\bar P$.
Denote by 
\begin{align*}
\bar T := T(\bar A,x^\ell,x^d).
\end{align*}
From Lemma~\ref{lem 1}, applied to $\bar A$, $x^\ell$, and $x^d$, we know that there exists a finite subset $\bar V$ of $\Z^n$ such that $\bar T = \cone \bar V,$ and for every $v \in \bar V$, we have $\norminf{v} \le \Delta$.
Since $x^\ell-x^d \in \bar T$, Caratheodory's theorem implies that there exist $m\le n$ vectors $v^1,\dots,v^m \in \bar V$ and $m$ positive scalars $\gamma_1,\dots, \gamma_m$ such that 
\begin{align} 
\label{will need this}
    x^\ell-x^d = x^\ell =\sum_{i=1}^{m}\gamma_{i}v^{i}.
\end{align}

The following simple observation will be used twice in our proof.
\begin{observation}
\label{ob1}
For $i=1,\dots,m,$ let $\lambda_i \in \R$,  such that $0\le \lambda_i \le \gamma_i$.
Then the vector $\sum_{i=1}^{m}\lambda_iv^i$ is in $\bar P.$
\end{observation}

\begin{prf}
Let $x:=\sum_{i=1}^{m}\lambda_iv^i$.
Since $x^d$ is the origin, we can write $x= x^d + \sum_{i=1}^{m}\lambda_iv^i$.
Using \eqref{will need this}, we can write $x$ also in the form
$
x=x^\ell-\sum_{i=1}^{m}(\gamma_i-\lambda_i) v^i.
$ 
From Lemma~\ref{lem 2}, applied with $\bar P$ and $\bar T$, we obtain $x \in \bar P.$
\end{prf}

We are now ready to define the vector $x^*$ as
\begin{align*}
    x^*:=\sum_{i=1}^{m}\floor{\gamma_i} v^{i}.
\end{align*}

\subsection{Properties of the vector $x^*$}
\label{sec property 2}

Note that $x^* \in \Z^n$ because $\floor{\gamma_i}$ and $v^{i}$ are all integer. 
From Observation~\ref{ob1}, we have $x^* \in \bar P$.

The next claim introduces several properties of $x^*$ that will be used later.

\begin{claim}
\label{claim xstar}
The vector $x^*$ satisfies the following properties:
\begin{enumerate}[label={(\alph*)}]
\setcounter{enumi}{3}
\item
\label{xstar d}
$\abs{x_i^*} \ge \chi_{\ell+1}-n\Delta$ for every $i \in N^\ell$;
\item
\label{xstar e}
$\{i \in \{1,\dots,k\} \mid x^*_i=0\} = Z^\ell;$
\item
\label{xstar f}
$\norminf{x^c-x^*} \le \psi_\ell+n\Delta$;
\item
\label{xstar g}
$x^c-x^* \in P$.
\end{enumerate}
\end{claim}

\begin{prf}
In this proof we will be using the upper bound on $\norminf{x^\ell-x^*}$ given by
\begin{align}
\label{eq cnc7}
\norminf{x^\ell-x^*}=\norminfL{\sum_{i=1}^{m}(\gamma_i-\floor{\gamma_i}) v^i} \le \sum_{i=1}^{m}\norminf{v^i}\le m\Delta\le n\Delta.
\end{align}
Next, we prove the properties of $x^*$ in the statement of the claim.

\smallskip

\ref{xstar d}. 
If $Z^\ell=\{1,\dots,k\}$ we are done, thus we assume $N^\ell \neq \emptyset.$
Let $i \in N^\ell$.
We have
\begin{align*}
\abs{x_i^\ell} = \abs{(x^\ell_i-x_i^*)+x^*_i} \le \abs{x^\ell_i-x_i^*} + \abs{x^*_i}.
\end{align*}
According to our assumption \ref{xell c-2}, we have $\abs{x^\ell_i} > \chi_{\ell+1}$, thus
\begin{align*}
    \abs{x_i^*} 
    \ge \abs{x^\ell_i}-\abs{x^\ell_i-x_i^*} 
    > \chi_{\ell+1}-\norminf{x^\ell-x^*}
    \ge \chi_{\ell+1}-n\Delta,
\end{align*}
where the last inequality holds by \eqref{eq cnc7}.

\smallskip

\ref{xstar e}.
Since the inequalities $x_i=0$, for $i \in Z^\ell,$ are valid for $\bar P$, and $x^* \in \bar P$, we know $\{i \in \{1,\dots,k\} \mid x^*_i=0\} \supseteq Z^\ell.$
One the other hand, given any index $i \in N^\ell$, we know from \ref{xstar d} that $\abs{x^*_i}\ge \chi_{\ell+1}-n\Delta>0$.
So $\{i \in \{1,\dots,k\} \mid x^*_i=0\} \subseteq Z^\ell.$
Thus we conclude that $\{i \in \{1,\dots,k\} \mid x^*_i=0\} = Z^\ell.$ 

\smallskip

\ref{xstar f}. This property follows directly from \ref{xell b} and \eqref{eq cnc7} as follows:
\begin{align*}
    \norminf{x^c-x^*}\le\norminf{x^c-x^\ell}+\norminf{x^\ell-x^*} \le \psi_\ell+n\Delta.
\end{align*}

\smallskip

\ref{xstar g}. 
Using the definition of $x^*$, the vector $x^c-x^*$ can be written as
\begin{align*}
    & x^c-x^*
    =x^c-\sum_{i=1}^m \floor{\gamma_i} v^i\\
    & x^c-x^*
    =(x^c-x^\ell)+(x^\ell-x^*)
    =(x^c-x^\ell)+\sum_{i=1}^m (\gamma_i - \floor{\gamma_i}) v^i.
\end{align*}
Recall that we have $\gamma_i > 0$ for every $i=1,\dots,m$.
Furthermore, we have $x^c \in P$ and, from \ref{xell a}, $x^c-x^\ell \in P$.
Let 
$T := T(A,x^\ell,x^d),$
and note that $\bar T \subseteq T,$ which implies $v^i \in T$ for every $i=1,\dots,m$.
Thus, from Lemma~\ref{lem 2}, applied with $P$ and $T$, we obtain that $x^c-x^* \in P.$
\end{prf}

Notice that \ref{xstar f} implies that $ \norminf{x^c-x^*} $ can be upper bounded by a function of $ n,\Delta,k,\epsilon $. In the next section, we will use \ref{xstar e} and \ref{xstar g} to show that the distance between $x^*$ and $x^d$ mainly depends on $ \abs{x^*_i} $ for $ i \in N^\ell $. In particular, when $ \abs{x^*_i} $ is large enough for every $ i \in N^\ell $, the vector $ x^* $ is a suitable approximation to $ x^d $. Together with \ref{xstar d}, this will imply that $ x^* $ is an $ \epsilon $-approximate solution to \eqref{pr IQP}.

\subsection{$x^*$ is an $\epsilon$-approximate solution to \eqref{pr IQP}}\label{sec analysis 1}

In this section we show that the vector $x^*$
is an $\epsilon$-approximate solution to \eqref{pr IQP}.
In Section~\ref{sec ub} we provide an upper bound for $ f(x^*)-f(x^d)$, while in Section~\ref{sec lb} we derive a lower bound for $ f_{\max}^d-f(x^d)$, where $ f_{\max}^d $ is the maximum value of $f$ on $P \cap \Z^n$.
In Section~\ref{sec ratio}, we use the two bounds to show that $x^*$ is an $\epsilon$-approximate solution to \eqref{pr IQP}.

\subsubsection{Upper bound on $ f(x^*)-f(x^d) $}
\label{sec ub}


\begin{claim}
\label{claim ratio 1}
We have $f(x^*)-f(x^d) \le 2(\psi_\ell+ n \Delta) \sum_{i \in N^\ell} q_i \abs{x^*_i}.$
\end{claim}

\begin{prf}
For ease of notation, let $\delta := x^c - x^*.$
We have 
\begin{align*}
    f(x^c) 
    & = f(\delta + x^*) 
    = \pare{ \sum_{i=1}^{k}-q_i \delta_i^2 + h^\transp \delta }
    + \pare{ \sum_{i=1}^{k}-q_i (x^*_i)^2+h^\transp x^* }
    - 2\sum_{i=1}^{k} q_i\delta_i x^*_i \\
    &=f(\delta)+f(x^*)- 2\sum_{i \in N^\ell} q_i\delta_i x^*_i 
    \ge f(x^c)+f(x^*)- 2\sum_{i \in N^\ell} q_i\delta_i x^*_i.
\end{align*}
In the last equality we used \ref{xstar e}.
    Furthermore, in the last inequality we used $f(\delta) \ge f(x^c)$ since $\delta \in P$ from \ref{xstar g} and $x^c$ is an optimal solution to \eqref{pr QP}.
	We obtain
	\begin{align*} 
	f(x^*) 
	\le 2 \sum_{i \in N^\ell} q_i \delta_i x^*_i \le 2 \sum_{i \in N^\ell} q_i \abs{\delta_i} \abs{x^*_i}  
	\le 2 (\psi_\ell+n\Delta) \sum_{i \in N^\ell} q_i \abs{x^*_i},
	\end{align*}
	where the last inequality holds from \ref{xstar f}.
	The claim follows by recalling that $f(x^d) = 0$ by Claim~\ref{claim eqt}.
\end{prf}

\subsubsection{Lower bound on $ f_{\max}^d-f(x^d) $}
\label{sec lb}

In this section we give a lower bound on $f_{\max}^d-f(x^d)$.
In our derivation, a fundamental role is played by the midpoint of $x^d$ and $x^*$, which we denote by $x^\triangle$, i.e., 
\begin{align*} 
x^\triangle 
:= \frac{x^d+x^*}{2}.
\end{align*} 
We first give a lower bound on $f(x^\triangle)-f(x^d)$.

\begin{claim}
\label{claim a}
We have $f(x^\triangle)-f(x^d)\ge \frac{1}{4}\sum_{i \in N^\ell} q_i(x^*_i)^2.$
\end{claim}
\begin{prf}
	The claim can be derived as follows:
	\begin{align*}
	f(x^\triangle) & = f \pare{ \frac{x^*}{2} }
	 = \frac 14 \sum_{i=1}^{k}-q_i (x^*_i)^2 + \frac 12 h^\transp x^* \\
	 & = \pare{ \frac{1}{2}\sum_{i=1}^{k}-q_i(x^*_i)^2+\frac{1}{2}h^\transp x^* } +\frac{1}{4}\sum_{i=1}^{k}q_i(x^*_i)^2 \\
	& = \frac{1}{2}f(x^*) + \frac{1}{4}\sum_{i \in N^\ell} q_i(x^*_i)^2 
	\ge \frac{1}{4}\sum_{i \in N^\ell} q_i (x^*_i)^2.
	\end{align*}
	In the last equality we used \ref{xstar e}.
	In the last inequality, we used $f(x^*) \ge f(x^d) = 0$, which holds because $x^d$ is optimal to \eqref{pr IQP} and $x^*$ is feasible to the same problem.
\end{prf}

Recall that the goal of this section is to obtain a lower bound on $f_{\max}^d-f(x^d)$.
Since both $x^d$ and $x^*$ are in $P$, the vector $x^\triangle$ is in $P$ as well.
Therefore, if $x^\triangle \in \Z^n$, then $f_{\max}^d \ge f(x^\triangle)$, and the bound of Claim~\ref{claim a} yields a bound on $f_{\max}^d-f(x^d)$.
However, $x^\triangle$ is not always an integer vector.
Thus we define two integer points $x^l$ and $x^r$ whose midpoint is $x^\triangle$: 
\begin{align*}
x^l & := 
x^d+\sum_{i \mid \floor{\gamma_i} \text{ odd}} \frac{\floor{\gamma_i}-1}{2} v^i
+\sum_{i \mid \floor{\gamma_i} \text{ even}} \frac{\floor{\gamma_i}}{2} v^i\\
x^r & := 
x^d+\sum_{i \mid \floor{\gamma_i} \text{ odd}} \frac{\floor{\gamma_i}+1}{2} v^i
+\sum_{i \mid \floor{\gamma_i} \text{ even}} \frac{\floor{\gamma_i}}{2} v^i.
\end{align*}

We now show that both $x^l$ and $x^r$ are in $\bar P \cap \Z^n$.
Clearly, $0 \le \frac{\floor{\gamma_i}}{2} \le \gamma_i$.
Furthermore, if $\floor{\gamma_i}$ is odd, we have $\floor{\gamma_i} \ge 1$, which implies $0 \le \frac{\floor{\gamma_i}-1}{2} \le \frac{\floor{\gamma_i}+1}{2} \le \gamma_i.$ 
By Observation~\ref{ob1}, we know that both $x^l$ and $x^r$ are in $\bar P$.
Since all coefficients $\frac{\floor{\gamma_i} \pm 1}{2}$ and  $\frac{\floor{\gamma_i}}{2}$ are integer, we conclude that both $x^l$ and $x^r$ are in $\bar P \cap \Z^n$.

Let $D \subset \R^n$ be the smallest box containing $x^l$ and $x^r$, i.e.,
\begin{align*}
D := [\min\{x^l_1,x^r_1\}, \max\{x^l_1,x^r_1\}] \times \dots \times [\min\{x^l_n,x^r_n\}, \max\{x^l_n,x^r_n\}].
\end{align*}

In the reminder of the proof we denote by $q : \R^n \to \R$ the quadratic part of the objective function $f$, i.e., 
\begin{align*}
q(x):=\sum_{i=1}^{k}-q_ix_i^2.
\end{align*} 
We also define the affine function $\lambda : \R^n \to \R$ which achieves the same value as $q$ at the vertices of the box $D$:
\begin{align*}
\lambda(x) := \sum_{i=1}^{k} -q_i \big( (x^l_i+x^r_i)x_i - x^l_ix^r_i \big).
\end{align*}
We have the following claim.

\begin{claim}
\label{claim lambdabound}
    For every $x\in D$ we have 
	$
	\lambda(x)\le q(x) 
	\le \lambda(x)+\frac{(n\Delta)^2}{4}\sum_{i \in N^\ell} q_i. 
	$
\end{claim}

\begin{prf}
	Since $\lambda$ achieves the same value as $q$ at each vertex of $D$ and $q$ is a concave function, we have 
	$\lambda(x) \le q(x)$, for every $x \in D$.

	Using the definitions of $q$ and $\lambda$ we obtain
	\begin{align*}
	q(x)-\lambda(x) & 
	= \sum_{i=1}^{k} -q_i \pare{ x_i^2 - (x^l_i+x^r_i)x_i + x^l_ix^r_i }
	=\sum_{i=1}^{k}-q_i(x_i-x^l_i)(x_i-x^r_i)  \\
	& \le \frac{1}{4}\sum_{i \in N^\ell} q_i(x^l_i-x^r_i)^2.
	\end{align*}
	The inequality holds because, for each $i=1,\dots,k$, the univariate quadratic function $-q_i(x_i-x^l_i)(x_i-x^r_i)$ achieves its maximum at $ \frac{x^l_i+x^r_i}{2}.$ 
	In particular, if $i \in Z^\ell,$ the maximum is 0.
	This is because both $x^l$ and $x^r$ are in $\bar P,$ which implies $x_i^l = x_i^r = 0$ for every $i \in Z^\ell$.
	
	From the definition of $x^l$ and $x^r$, we obtain $\norminf{x^r-x^l} \le \norminf{\sum_{i=1}^{m} v^i} \le m\Delta \le n\Delta$.
	Therefore, we have $q(x)\le \lambda(x)+ \frac{(n\Delta)^2}{4} \sum_{i \in N^\ell} q_i$.
\end{prf}

\begin{claim}
\label{claim b}
There exists $\tilde x \in \{x^l, x^r\}$ such that 
$
f(\tilde{x}) - f(x^\triangle) \ge  - \frac{(n\Delta)^2}{4}\sum_{i \in N^\ell} q_i.
$
\end{claim}

\begin{prfc}
Let $g : \R^n \to \R$ be defined by $g(x):=\lambda(x)+h^\transp x$. 
Claim~\ref{claim lambdabound} implies that, for every $x\in D$, we have
\begin{equation*}
g(x)\le f(x) \le g(x)+\frac{(n\Delta)^2}{4}\sum_{i \in N^\ell}q_i.
\end{equation*}

Since $g$ is a linear function and $x^\triangle$ is the midpoint of $x^l$ and $x^r$, we know that $g(x^\triangle) \le g(\tilde x)$ for some $\tilde x \in \{ x^l,  x^r \}$. 
We derive the following relation: 
\begin{equation*}
f(x^\triangle) \le g(x^\triangle)+\frac{(n\Delta)^2}{4}\sum_{i \in N^\ell}q_i 
\le g(\tilde{x})+\frac{(n\Delta)^2}{4}\sum_{i \in N^\ell}q_i
\le f(\tilde{x})+\frac{(n\Delta)^2}{4}\sum_{i \in N^\ell}q_i. \qedhere
\end{equation*}
\end{prfc}

We are finally ready to state our lower bound on $f_{\max}^d-f(x^d)$.
\begin{claim}
\label{claim ratio 2}
	We have
$
	f_{\max}^d-f(x^d)\ge \frac{1}{4} \sum_{i \in N^\ell} q_i \big( (x^*_i)^2-(n\Delta)^2 \big).
$
\end{claim}

\begin{prfc}
	Combining Claim~\ref{claim a} and Claim~\ref{claim b}, we have
	\begin{align*}
	f_{\max}^d - f(x^d) 
	& \ge f(\tilde{x}) - f(x^d) 
	= \pare{ f(\tilde{x})-f(x^\triangle) } + \pare{ f(x^\triangle)-f(x^d) } \\
	& \ge \frac{1}{4} \sum_{i \in N^\ell}
	q_i \pare{ (x^*_i)^2 - (n\Delta)^2 }. \qedhere
	\end{align*}
\end{prfc}

\subsubsection{$x^*$ is an $\epsilon$-approximate solution}
\label{sec ratio}

In order to prove that $x^*$ is an $\epsilon$-approximate solution, we 
first prove the following observation.

\begin{observation}\label{ob2}
Let $a_i,b_i>0$, for $i=1,\dots,k.$ 
Then $ \frac{\sum_{i=1}^{k}a_i}{\sum_{i=1}^{k}b_i}\le \max_{i = 1,\dots,k}\frac{a_i}{b_i} $.
\end{observation}
\begin{prf}
To prove this statement, let $j \in \{1,\dots,k\}$ such that $\frac{a_j}{b_j} = \max_{i = 1,\dots,k}\frac{a_i}{b_i}$.
	Then,
	\begin{align*}
	\frac{\sum_{i=1}^{k}a_i}{\sum_{i=1}^{k}b_i}-\frac{a_j}{b_j}
	=\frac{\sum_{i=1}^{k}(b_j a_i-a_j b_i)}{b_j \sum_{i=1}^{k}b_i}.
	\end{align*}
We only need to show that the right-hand side of the latter equation is nonpositive.
To see this, notice that,
	for $i=1,\dots,k$, we have $\frac{a_i}{b_i} \le \frac{a_j}{b_j}$, thus
	$b_j a_i - a_j b_i\le 0.$
\end{prf}

\begin{claim}
\label{claim solution}
The vector $x^*$ is an $\epsilon$-approximate solution to \eqref{pr IQP}.
\end{claim}

\begin{prf}
Consider first the case $Z^\ell=\{1,\dots,k\}.$ 
Then by \ref{xstar e}, we know that $x^*_i=0$, for $i=1,\dots,k.$ 
In this case, from Claim~\ref{claim ratio 1}, we know that $f(x^*)\le 0.$ 
By Claim~\ref{claim eqt}, this implies that $x^*$ is an optimal solution to \eqref{pr IQP}.

Now assume that $Z^\ell \subset \{1,\dots,k\}$, i.e., $N^\ell \neq \emptyset$. 
Observe that the quantity $f_{\max}^d-f(x^d)$ in the definition of $\epsilon$-approximate solution is positive.
This follows from Claim~\ref{claim ratio 2}, since for $i \in N^\ell$, we have $q_i > 0$ by assumption and $\abs{x_i^*} > n \Delta$ from \ref{xstar d}.
Therefore, we consider the ratio $\frac{f({x^*}) - f(x^d)}{f_{\max}^d-f(x^d)}$, and our aim is to show that it is upper bounded by $\epsilon$.
Using Claim~\ref{claim ratio 1}, Claim~\ref{claim ratio 2}, and Observation~\ref{ob2}, we derive the following bound: 
\begin{align*}
\frac{f({x^*}) - f(x^d)}{f_{\max}^d-f(x^d)} 
& \le 8(\psi_\ell+ n \Delta) \frac{\sum_{i \in N^\ell}q_i\abs{x_i^*}}{\sum_{i \in N^\ell}q_i((x^*_i)^2-(n\Delta)^2 )} \\
& \le 8(\psi_\ell+ n \Delta) \max_{i \in N^\ell} \frac{\cancel{q_i}\abs{x_i^*}}{\cancel{q_i}((x^*_i)^2-(n\Delta)^2)}.
\end{align*}
In particular, the latter max can be written in the form
\begin{align*}
\max_{i \in N^\ell}\frac{\abs{x_i^*}}{(x^*_i)^2-(n\Delta)^2}
= \max_{i \in N^\ell} \frac{1}{\abs{x^*_i} - \frac{(n\Delta)^2}{\abs{x^*_i}}}.
\end{align*}
In the right-hand side, the denominator is always positive due to $\abs{x_i^*} > n \Delta$.
Let $s$ be the index in $\{1,\dots,k\}$ that achieves $\min_{i \in N^\ell}\abs{x_i^*}$.
Note that this index exists because of our assumption $N^\ell \neq \emptyset$.
Then the max is achieved by the index $s$.
In fact, in the denominator in the right-hand side, the term $\abs{x^*_i}$ is minimized by $s$, while the term $\frac{(n\Delta)^2}{\abs{x^*_i}}$ is maximized by $s$.

From \ref{xstar d}, we have 
$\abs{x_s^*} \ge \chi_{\ell+1}-n\Delta= \frac{8}{\epsilon}(\psi_\ell+n\Delta) + n\Delta,$ where the equality can be obtained using the definition of $\chi_{\ell+1}$ and of $\psi_\ell$.
We obtain
\begin{align*}
\frac{f({x^*}) - f(x^d)}{f_{\max}^d-f(x^d)} 
& \le 
\frac{8(\psi_\ell+ n \Delta)}{\abs{x^*_s} - \frac{(n\Delta)^2}{\abs{x^*_s}}}
\le \frac{8(\psi_\ell+ n \Delta)}{\frac{8}{\epsilon}(\psi_\ell+ n \Delta)+n\Delta-\frac{(n\Delta)^2}{\frac{8}{\epsilon}(\psi_\ell+ n \Delta)+n\Delta} } \\
&
=\frac{\cancel{8(\psi_\ell+n\Delta)}\pare{\frac{8}{\epsilon}(\psi_\ell+n\Delta)+n\Delta}}{\cancel{8(\psi_\ell+n\Delta)}\pare{\frac{8}{\epsilon^2}(\psi_\ell+n\Delta)+\frac{2}{\epsilon}n\Delta}}  \\
& =\epsilon\frac{8(\psi_\ell+n\Delta)+n\Delta\epsilon}{8(\psi_\ell+n\Delta)+2n\Delta\epsilon} < \epsilon,
\end{align*}
where the first equality can be obtained by multiplying the numerator and the denominator by $\frac{8}{\epsilon}(\psi_\ell+ n \Delta)+n\Delta.$
This implies that $x^*$ is an $\epsilon$-approximate solution to \eqref{pr IQP}.
\end{prf}

From \ref{xstar f}, we know that
$\norminf{x^c-x^*}\le \psi_\ell+n\Delta \le \psi_k+n\Delta.$
Furthermore, from Claim~\ref{bound}, we have $\psi_k + n\Delta \le n\Delta(\frac{10\Delta}{\epsilon}+1)^k.$
This completes the proof of Theorem~\ref{th main}\ref{th main 1}.
Therefore, in the remainder of the proof we only need to show Theorem~\ref{th main}\ref{th main 2}.

	\subsection{Construction of the vector $x^\star$}
		\label{sec star}

In this section we introduce the vector $x^\star$ in the statement of Theorem~\ref{th main}\ref{th main 2}.
The point $x^\star$ is defined by
	\begin{align*}
	    x^\star:=x^c-x^*.
	\end{align*}

	From the definition of $x^*$, we obtain 
	\begin{align*}
	x^\star
	=x^c-\sum_{i=1}^{m}\floor{\gamma_i} v^i=
	x^c-x^\ell+\sum_{i=1}^{m}(\gamma_i-\floor{\gamma_i})v^i.
	\end{align*}
	We know that $x^c \in P,$ and, from \ref{xell a}, we know that $x^c-x^\ell \in P$ as well.
	Thus, from Lemma~\ref{lem 2}, applied with $P$ and $T,$ we know that $x^\star \in P.$

	\subsection{$x^\star$ is an $\epsilon$-approximate solution to \eqref{pr QP}}\label{sec analysis 2}

	In this section we show that the vector $x^\star$ is an $\epsilon$-approximate solution to \eqref{pr QP}. 
	To do this, we first give an upper bound on $f(x^\star)-f(x^c)$, and then a lower bound on $f_{\max}^c-f(x^c)$, where $f_{\max}^c$ is the maximum value of $f$ on $P$.
	The two bounds are then used to show that $x^\star$ is an $\epsilon$-approximate solution to \eqref{pr QP}.
	
	\begin{claim}\label{rub}
We have $f(x^\star)-f(x^c) \le 2(\psi_\ell+ n \Delta) \sum_{i \in N^\ell} q_i \abs{x^*_i}.
	$
	\end{claim}

	\begin{prf}
	First, we derive an upper bound on $f(x^\star)$.
	According to the definition of $x^\star$, we get
	\begin{align*}
	    f(x^\star) & = f(x^c-x^*) \\
& 	                 = \pare{-\sum_{i=1}^{k}q_i(x^c_i)^2+h^\transp x^c} + \pare{-\sum_{i=1}^{k}q_i(x^*_i)^2-h^\transp x^*}+
	                 2\sum_{i=1}^{k}q_ix^c_ix^*_i \\
	                 & =f(x^c)+f(-x^*)+2\sum_{i \in N^\ell}q_ix^c_ix^*_i \\
	                 & =f(x^c)+f(-x^*)+2\sum_{i \in N^\ell}q_i(x^*_i+x^\star_i)x^*_i,
	\end{align*}
	where in the third equality we used \ref{xstar e}.
	
	To derive from the above formula an upper bound on $f(x^\star)-f(x^c)$, we need to upper bound $f(-x^*)$.
	Since $x^d$ is the optimal solution to \eqref{pr IQP} and $x^* \in P\cap \Z^n,$ we know that
	\begin{align*}
	f(x^*)=-\sum_{i=1}^k q_i(x^*_i)^2+h^\transp x^* \ge f(x^d)=0. 
	\end{align*}
	Thus, we get
	\begin{align*}
	f(-x^*)
	& =-\sum_{i=1}^{k}q_i(x^*_i)^2-h^\transp x^*
	\le f(x^*) - \sum_{i=1}^{k}q_i(x^*_i)^2-h^\transp x^*\\
	& =-2\sum_{i=1}^{k}q_i(x^*_i)^2
	=-2\sum_{i \in N^\ell}q_i(x^*_i)^2.
	\end{align*}
	
	We obtain
	\begin{align*}
	f(x^\star)-f(x^c) & 
	\le -2\sum_{i \in N^\ell}q_i(x^*_i)^2+2\sum_{i \in N^\ell}q_i(x^*_i+x^\star_i)x^*_i
	=2\sum_{i \in N^\ell}q_ix^\star_ix^*_i \\
	& \le 2\sum_{i \in N^\ell}q_i \abs{x^\star_i} \abs{x^*_i}
	\le 2(\psi_\ell+ n \Delta) \sum_{i \in N^\ell} q_i \abs{x^*_i}. 
	\end{align*}
	The last inequality holds because, from \ref{xstar f}, we have $\norminf{x^\star}=\norminf{x^c-x^*}\le \psi_\ell+n\Delta.$
	\end{prf}

\begin{claim}\label{rlb}
We have
$
f_{\max}^c-f(x^c)\ge \frac{1}{4}\sum_{i \in N^\ell}q_i(x^*_i)^2.
$
\end{claim}

\begin{prfc}
Define the midpoint of $x^c$ and $x^\star$ as
\begin{align*}
x^\diamond:=\frac{x^c+x^\star}{2}.
\end{align*}
Then
\begin{align*}
	f(x^\diamond) & = f \pare{ \frac{x^c+x^\star}{2} }
	= \frac{1}{4}\sum_{i=1}^{k}-q_i \big( (x^c_i)^2+2x^c_ix^\star_i+(x^\star_i)^2 \big) +\frac 12 h^\transp ( x^c + x^\star) \\
	& = \pare{ \frac{1}{2}\sum_{i=1}^{k}-q_i(x^c_i)^2+\frac{1}{2}h^\transp x^c } + \pare{ \frac{1}{2}\sum_{i=1}^{k}-q_i(x^\star_i)^2+\frac{1}{2}h^\transp x^\star } +\frac{1}{4}\sum_{i=1}^{k}q_i(x^\star_i-x^c_i)^2 \\
	& = \frac{1}{2}f(x^c)+\frac{1}{2}f(x^\star) + \frac{1}{4}\sum_{i \in N^\ell}q_i(x^*_i)^2 
	 \ge f(x^c) + \frac{1}{4}\sum_{i \in N^\ell} q_i (x^*_i)^2.
	\end{align*}
In the last inequality, we used $f(x^\star) \ge f(x^c)$, which holds because $x^c$ is optimal to \eqref{pr QP} and $x^\star$ is feasible to the same problem.
Since $x^\diamond$ is the midpoint of $x^\star$ and $x^c$, we know that $x^\diamond \in P.$ 
Thus we have
\begin{equation*}
f_{\max}^c-f(x^c)\ge f(x^\diamond)-f(x^c) \ge \frac{1}{4}\sum_{i \in N^\ell} q_i (x^*_i)^2. \qedhere
\end{equation*}
\end{prfc}

	\begin{claim}\label{rresult}
The vector $x^\star$ is an $\epsilon$-approximate solution to \eqref{pr QP}.
	\end{claim}
	
	\begin{prf}
	As in the proof of Claim~\ref{claim solution}, it is simple to check that the quantity $f_{\max}^c-f(x^c)$ in the definition of $\epsilon$-approximate solution is positive.
	This allows us to consider the ratio $\frac{f({x^\star}) - f(x^c)}{f_{\max}^c-f(x^c)}$, and our aim is to show that it is upper bounded by $\epsilon$.

	Using Claim~\ref{rub}, Claim~\ref{rlb}, and Observation~\ref{lem 2}, we can derive the following bound: 
\begin{align*}
\frac{f({x^\star}) - f(x^c)}{f_{\max}^c-f(x^c)} 
& \le 8(\psi_\ell+ n \Delta) \frac{\sum_{i \in N^\ell}q_i\abs{x_i^*}}{\sum_{i \in N^\ell}q_i(x^*_i)^2} \\
& \le 8(\psi_\ell+ n \Delta) \max_{i \in N^\ell} \frac{\cancel{q_i}\abs{x_i^*}}{\cancel{q_i}(x^*_i)^2}=\frac{8(\psi_\ell+ n \Delta)}{\abs{x_s^*}},
\end{align*}
where $s$ is an index such that $\abs{x^*_s} = \min \{\abs{x^*_i} \mid i \in N^\ell\}.$

From \ref{xstar d}, we have 
$\abs{x_s^*} \ge \chi_{\ell+1}-n\Delta$, and using the definition of $\chi_{\ell+1}$ the latter quantity equals $n\Delta +\frac{8}{\epsilon}(\psi_\ell+n\Delta).$ 
We get
\begin{align*}
    \frac{8(\psi_\ell+ n \Delta)}{\abs{x^*_s}}
    \le \frac{8(\psi_\ell+ n \Delta)}{n\Delta +\frac{8}{\epsilon}(\psi_\ell+n\Delta)}<\frac{8(\psi_\ell+ n \Delta)}{\frac{8}{\epsilon}(\psi_\ell+n\Delta)}= \epsilon.
\end{align*}
This implies that $x^\star$ is an $\epsilon$-approximate solution to \eqref{pr QP}.
\end{prf}

  From the definition of $x^\star$ and from \ref{xstar f}, we obtain
    \begin{align*}
    \norminf{x^\star-x^d}= \norminf{x^c-x^*}\le \psi_\ell+n\Delta \le \psi_k+n\Delta.
    \end{align*}
    Moreover, from Claim~\ref{bound}, we have $\psi_k+n\Delta \le n\Delta(\frac{10\Delta}{\epsilon}+1)^k.$ 
    This completes the proof of Theorem~\ref{th main}\ref{th main 2}, and of Theorem~\ref{th main}. 

\section{Lower bounds on the distance of solutions}
\label{sec tight}

In this section we discuss how far from optimal are the proximity bounds in Theorem~\ref{th main}.
The main ingredient in the derivation of our lower bounds is a polyhedron $\bar P$ that we introduce next.

\begin{definition}
\label{def poly}
For every 
$n, \Delta, t \in \Z$ with
$n\ge 1$, 
$\Delta \ge 1$,
$t \ge 0$, 
and
$\beta \in (0, 1)$,
let $\bar P \subset \R^n$ be the polyhedron defined by the following inequalities:
\begin{align*}
    & -t \le x_1-\Delta \sum_{i=2}^n x_i \le t \\
    & 0 \le x_i \le \beta & i=2,\dots, n.
\end{align*}
\end{definition}

Clearly, the polyhedron $\bar P$ has dimension $n$ if $t \ge 1$.
Note that $\bar P$ can be obtained from the polytope
\begin{align*}
\bar P_y := \{(y_1, x_2, \dots,x_n) \mid -t \le y_1 \le t, \ 0 \le x_i \le \beta, i=2,\dots,n \}
\end{align*}
by replacing variable $y_1$ with $x_1 := y_1 + \Delta \sum_{i=2}^n x_i.$
The vertices of the polytope $\bar P_y$ are all 
vectors with components $y_1 = \pm t$, and $x_i \in \{0,\beta\}$ for $i=2,\dots,n$.
Therefore, $\bar P$ is bounded and its vertices are all 
vectors with components $x_i \in \{0,\beta\}$ for $i=2,\dots,n$, and component $x_1 = \pm t + \Delta \sum_{i=2}^n x_i$.
In particular, the vertex with the largest $x_1$ is
\begin{align*}
v:=\pare{t+(n-1)\beta\Delta,\beta,\dots,\beta},
\end{align*}
and will play an important role in our arguments.

Next, we focus on the integer points in $\bar P$.
Since $\beta < 1$, any vector in $\bar P \cap \Z^n$ satisfies $x_i=0$, for $i=2,\dots, n$, and the first component of such vectors are $-t,-t+1,\dots,t$.
Since $t \ge 0$, the set $\bar P \cap \Z^n$ contains the origin and is therefore nonempty. 
In particular, the integer point in $\bar P$ with the largest $x_1$ is
\begin{align*}
u:=(t,0,\dots,0).
\end{align*}
The vectors $u$ and $-u$ will often be used in the later proofs.


\begin{observation}
\label{obs delta}
Let $A$ be the constraint matrix defining $\bar P$.
Then each 
subdeterminant of $A$ is in $\{0, \pm 1, \pm \Delta\}$.
\end{observation}

\begin{prf}
The constraint matrix of the system defining $\bar P$ is
\begin{align*}
A = 
\begin{pmatrix}
1 & -\Delta_{n-1}^\transp \\
-1 & \Delta_{n-1}^\transp \\
0_{n-1} & I_{n-1} \\
0_{n-1} & -I_{n-1}
\end{pmatrix},
\end{align*}
where $I_{n-1}$ denotes the $(n-1) \times (n-1)$ identity matrix, and $0_{n-1}$ (resp.~$\Delta_{n-1}$) denotes the $(n-1)$-dimensional vector with all entries equal to zero (resp.~$\Delta$).
Let $d$ be the determinant of a square submatrix $M$ of $A$.
If $M$ has linearly dependent rows, then $d=0$.
Thus we now assume that $M$ does not have linearly dependent rows.
Up to mutiplying rows of $M$ by $-1$, which is an operation that can only change the sign of the determinant $d$, the matrix $M$ is a submatrix of 
\begin{align*}
\begin{pmatrix}
1 & -\Delta_{n-1}^\transp \\
0_{n-1} & I_{n-1} \\
\end{pmatrix}.
\end{align*}
It is well known that adding unit rows to a matrix can only add $0,1$ and change the sign to its possible 
subdeterminants.
Therefore, either $d \in \{0,1\}$, or $d$ is a 
subdeterminant
of the matrix 
$
\begin{pmatrix}
1 & -\Delta_{n-1}^\transp
\end{pmatrix}.
$
The latter matrix has only one row and its 
subdeterminants
are $1, -\Delta$.
\end{prf}

For brevity, in this section, we say that a \eqref{pr IQP} or \eqref{pr QP} has \emph{subdeterminant $\Delta$} if the maximum of the absolute values of the subdeterminants of the constraint matrix $A$ is $\Delta$.

\subsection{Tightness in Integer Linear Programming}

In this section we consider our problems \eqref{pr IQP} and \eqref{pr QP} under the additional assumption $k=0$.
In this special case, \eqref{pr IQP} is a general Integer Linear Programming (ILP) problem, while \eqref{pr QP} is the corresponding Linear Programming (LP) problem, also known as the standard linear relaxation of (ILP).

We remark that, for $k=0$, Theorem~\ref{th main} reduces to the proximity bound by Cook et al.~\cite{CooGerSchTar86} for Integer Linear Programming. 
In particular, this result yields the upper bound
\begin{align*}
\min \bra{ \norminf{x^c-x^d} \mid x^d \text{ opt.~to (ILP)}, \ x^c \text{ opt.~to (LP)} } \le n \Delta.
\end{align*}

The impact of the polytope $\bar P$ 
is immediately apparent, as it allows us to prove that the above upper bound
$n \Delta$ is asymptotically best possible.
To the best of our knowledge this tightness result was previously known only for $\Delta = 1$ \cite{SchBookIP,PaaWeiWel18}.

\begin{proposition}
\label{prop ILP ub}
For every 
$n, \Delta \in \Z$ with
$n\ge 1$, 
$\Delta \ge 1$, and
$\beta \in (0, 1)$,
there exists an instance of (ILP) 
with subdeterminant $\Delta$ for which
\begin{align*}
\min \bra{ \norminf{x^c-x^d} \mid x^d \text{ opt.~to \textnormal{(ILP)}}, \ x^c \text{ opt.~to \textnormal{(LP)}}} = (n-1)\beta\Delta \in \Omega(n\Delta).
\end{align*}
\end{proposition} 

\begin{prf}
Let $n$, $\Delta$, $\beta$ be as in the statement.
Consider the (ILP) problem
\begin{align}
\label{eq IQP tight}
\begin{split}
	\max \ & x_1 \\
	\stt \ & x \in \bar P \cap \Z^n,
\end{split}
\end{align}
where the parameter $t$ in the definition of $\bar P$ can be chosen to be any integer greater than or equal to zero.
From Observation~\ref{obs delta}, problem \eqref{eq IQP tight} has subdeterminant $\Delta$.
The unique optimal solution of \eqref{eq IQP tight} is the vector $u$, while the unique optimal solution of the corresponding (LP) is the vertex $v$ of $\bar P$.
We obtain $\norminf{v-u} = (n-1)\beta\Delta$.
\end{prf}

\subsection{
Lower bounds in Integer Quadratic Programming}

Let us now get back to the general case of \eqref{pr IQP} where $k$ can be positive. 
In this setting, even for $k=1$, problem \eqref{pr ex} in Example~\ref{ex no bound} shows that it is not possible to upper bound the distance
\begin{align*}
\min \bra{ \norminf{x^c-x^d} \mid x^d \text{ opt.~to } \eqref{pr IQP}, \ x^c \text{ opt.~to } \eqref{pr QP} }
\end{align*}
with a function that depends only on $n$ and $\Delta$.
Therefore, we focus instead on the two quantities 
\begin{align*}
\delta^*_\epsilon & 
:= \min \bra{ \norminf{x^c-x^*} \mid x^* \text{ $\epsilon$-approx.~to } \eqref{pr IQP}, \ x^c \text{ opt.~to } \eqref{pr QP} }, \\
\delta^\star_\epsilon & := 
\min \bra{ \norminf{x^\star - x^d} \mid  x^d \text{ opt.~to } \eqref{pr IQP}, \ x^\star \text{ $\epsilon$-approx.~to } \eqref{pr QP} }.
\end{align*}
Our Theorem~\ref{th main} implies that both $\delta^*_\epsilon$ and $\delta^\star_\epsilon$ are upper bounded by
\begin{align*}
n\Delta \pare{\frac{10\Delta}{\epsilon}+1}^k 
\quad \in O\pare{\frac{n \Delta^{k+1}}{\epsilon^k}}.
\end{align*}
In the next two sections we gain insight on how far from optimal are our proximity results.
This is done by providing lower bounds on both $\delta^*_\epsilon$, in Section~\ref{sec tight a}, and on $\delta^\star_\epsilon$, in Section~\ref{sec tight b}.
Note that our bounds can be further improved, as we are only interested here in the asymptotic behaviour of $\delta^*_\epsilon$ and $\delta^\star_\epsilon$.

We remark that the problems that we present in the following results are of the form \eqref{pr IQP} and \eqref{pr QP} with an additional constant in the objective function.
We decided to keep these constants to simplify the presentation, and we observe that the presence of these constants does not affect optimal or $\epsilon$-approximate solutions thanks to Lemma~\ref{lem tans}.

\subsubsection{Lower bounds on $\delta^*_\epsilon$}
\label{sec tight a}

To begin with, we present a special  \eqref{pr IQP} problem, which will be useful in the subsequent discussion.
For every $n, \Delta, t \in \Z$ with $n\ge 1$, $\Delta \ge 1$, $t \ge 0$, and $a \in \R$, $\beta \in (0, 1)$, consider the \eqref{pr IQP}
\begin{align}
\label{pr tight}
\begin{split}
	\min \ & f(x)=
	 -(x_1-a)^2
	- \frac{(t+n\Delta)^2}{\beta^2} \sum_{i=2}^n x_i^2  \\
	\stt \ & x \in \bar P \cap \Z^n,
	\end{split}
\end{align}
where the polytope $\bar P$ is given in Definition~\ref{def poly}.
Note that problem \eqref{pr tight} has $k=n$ and, from Observation~\ref{obs delta}, subdeterminant $\Delta$.

The next lemma provides some information about \eqref{pr tight} and its corresponding \eqref{pr QP}.
We remind the reader that the vectors $u,v$ are defined right after Definition~\ref{def poly}.

\begin{lemma}
\label{optial QP}
If $0 < a < (n-1)\beta\Delta$, then the vector $-u$ is the unique optimal solution to \eqref{pr tight}, and the vector $v$ is the unique optimal solution to the corresponding \eqref{pr QP}.
\end{lemma}

\begin{prf}
Consider problem \eqref{pr tight} and assume $0 < a < (n-1)\beta\Delta$.

We first show that the vector $-u$ is the unique optimal solution to \eqref{pr tight}.
We have seen that any vector in $\bar P \cap \Z^n$ satisfies $x_i = 0$, for $i=2,\dots,n$, and the first component ranges in $-t,-t+1,\dots,t$.
Our assumption $a > 0$ then implies that the vector $-u = (-t,0,\dots,0)$ is the unique optimal solution to \eqref{pr tight}. 

Next, we show that the vertex $v$ of $\bar P$ is the unique optimal solution to the corresponding \eqref{pr QP}.
Note that 
\begin{align*}
    f(v)& =-(t+(n-1)\beta \Delta -a)^2-(n-1)(t+n\Delta)^2.
\end{align*}
Since $\bar P$ is a polytope and the objective is concave, we only need to show that any other vertex $v'$ of $\bar P$ has cost strictly larger than $v$.

First, assume that $v'_i = \beta$ for every $i=2,\dots,n$.
Then $v'_1 = -t+(n-1)\beta \Delta.$ 
Notice that, if $t=0$, then $v'$ and $v$ are the same point.
Therefore, we assume that $t > 0$.
We have
\begin{align*}
    f(v')& =-(-t+(n-1)\beta \Delta -a)^2-(n-1)(t+n\Delta)^2.
\end{align*}
Since $(n-1)\beta \Delta -a > 0$, due to the fact that $t > 0$, we obtain $\abs{t+(n-1)\beta \Delta -a} > \abs{-t+(n-1)\beta \Delta -a}$.
We have therefore shown $f(v) < f(v')$.

We can now assume that $v'_i = 0$ for some $i=2,\dots,n$.
If we denote by $m$ the number of components among $v'_2,\dots,v'_n$ that are equal to $\beta$, then we have $m \le n-2$.
Hence
\begin{align*}
    f(v')& =-(v'_1 - a)^2-m(t+n\Delta)^2 
    \ge -(v'_1 - a)^2-(n-2)(t+n\Delta)^2 \\
    & > -(t+n\Delta)^2-(n-2)(t+n\Delta)^2 
     = - (n-1)(t+n\Delta)^2 \ge f(v).
\end{align*}
Now we explain how we obtain the second inequality.
From the definition of $\bar P,$ we know that $-t \le v'_1 \le t+ (n-1)\beta \Delta.$ 
Since $0\le a \le (n-1)\beta\Delta,$ we get 
\begin{align*}
-t
-(n-1)\beta\Delta
\le v'_1 - a \le t+(n-1)\beta\Delta.
\end{align*} 
We obtain 
\begin{align*}
(v'_1 - a)^2 \le (t+(n-1)\beta\Delta)^2 < (t+n\Delta)^2,
\end{align*}
which implies the second inequality.
In this second case we have shown that $f(v) < f(v')$ holds for every $t \ge 0$.
This concludes the proof that $v$ is the unique optimal solution to the \eqref{pr QP} corresponding to \eqref{pr tight}.
\end{prf}

In the next proposition we highlight a key difference between Separable Concave Integer Quadratic Programming and Integer Linear Programming.
More in detail, we discuss an important difference between problem \eqref{pr IQP} with $k \ge 1$ and the same problem with $k=0$.
Consider a feasible instance of \eqref{pr IQP}, and let $x^c$ be an optimal solution to the corresponding \eqref{pr QP}.
According to Cook et al.~\cite{CooGerSchTar86}, we can always find integer points $x$ in $P$ with $\norminf{x^c - x} \le n\Delta.$ 
Furthermore, if $k=0$, one of these vectors is optimal to \eqref{pr IQP}.
However, this is not true for the case $k \ge 1$. 
In fact, when 
$k \ge 1$, 
the set $\{x \in P \cap \Z^n \mid \norminf{x^c-x} \le n\Delta\}$ not only might contain no optimal solution to \eqref{pr IQP}, but it also might contain only arbitrarily bad solutions, i.e., vectors that are not $\epsilon$-approximate solution to \eqref{pr IQP}, for any $\epsilon \in (0,1)$.

\begin{proposition}
\label{prop IQP ub first}
For every $\Delta \in \Z$ with $\Delta \ge 1$, and 
$\epsilon \in (0,1),$ there is an instance of \eqref{pr IQP} 
with subdeterminant $\Delta$ 
and $k=n$
for which $\delta^*_\epsilon > n\Delta.$
\end{proposition}


\begin{prf}
Let $\Delta, \epsilon$ be as in the statement, and consider problem \eqref{pr tight} with $n := \ceil{\frac{4-3\sqrt{\epsilon}}{1-\sqrt{\epsilon}}} \ge 5$, $\beta := \frac{n-4}{n-3} \in [\frac 12 ,1)$, $a := (n-3)\beta\Delta$, and $t := \frac{n+(n-3)\beta}{2} \Delta \ge 3$. 
Note that, both $a$ and $t$ are integer.
From Lemma~\ref{optial QP}, we know that $x^d = -u$ is the unique optimal solution to \eqref{pr tight} and $x^c = v$ is the unique optimal solution to the corresponding \eqref{pr QP}. 

Let
$S:=\{x \in \bar P \cap \Z^n \mid \norminf{x^c-x} \le n\Delta\}$.
It suffices to show that there is no $\epsilon$-approximate solution to \eqref{pr tight} in $S$.
Using the definition of $v$, we derive
\begin{align*}
S =\{x \in \R^n \mid t-(n-(n-1)\beta)\Delta \le x_1 \le t, \ x_1 \in \Z, \ x_i=0, i=2,\dots,n \},
\end{align*}
and it can be checked that the quantity $t-(n-(n-1)\beta) \Delta$ is in $(-t,t)$.
Using $n \ge 5,$ it can be checked that $a$ is smaller than the midpoint between the two points $t-(n-(n-1)\beta)\Delta$ and $t$.
Due to the concavity of the objective, and the fact that $t \in \Z$, this implies that the vector $u$ is a minimizer of the objective function $f(x)$ over the set $S$.
Therefore, it suffices to show that the vector $u$ is not an $\epsilon$-approximate solution to \eqref{pr tight}.

Let $u':=(a,0,\dots,0).$ 
Since $-t<a<t,$ we have $u' \in \bar P.$ 
Moreover, since $a$ is integer, it is simple to check that $f^d_{\max}=0$ and it is achieved at $u'.$ 
Furthermore, we have 
\begin{align*}
f(u) & = -(t-a)^2 = -(t-(n-3)\beta\Delta)^2
=-\pare{\frac{n-(n-3)\beta}{2}\Delta}^2
=-4\Delta^2, \\
f(x^d) & = -(t+a)^2 =-(t+(n-3)\beta\Delta)^2
=-\pare{\frac{n+3(n-3)\beta}{2}\Delta}^2
=-(2n-6)^2\Delta^2.
\end{align*}
We obtain
\begin{align*}
\frac{f(u)-f(x^d)}{f^d_{\max}-f(x^d)}
=\frac{\cancel 4 (n-4)(n-2)\cancel{\Delta^2}}{\cancel 4 (n-3)^2\cancel{\Delta^2}} 
>\frac{(n-4)^2}{(n-3)^2}
=\pare{1-\frac{1}{n-3}}^2 \ge \epsilon,
\end{align*}
where the last inequality can be checked by plugging in $n.$
Therefore, the vector $u$ is not an $\epsilon$-approximate solution to \eqref{pr tight}.
\end{prf}


In particular, Proposition~\ref{prop IQP ub first} 
shows that if $k=n$, then $\delta^*_\epsilon$ can grow at least linearly with respect to both $n$ and $\Delta$.
In the next proposition, we use Lemma~\ref{optial QP} to derive our main lower bound on $\delta^*_\epsilon$.

\begin{proposition}
\label{prop IQP ub}
For every $n, \Delta \in \Z$ with $n\ge 2$, $\Delta \ge 1$, and $\epsilon \in (0,1]$, 
there exists an instance of \eqref{pr IQP} 
with subdeterminant $\Delta$ 
and $k=n$
for which
\begin{align*}
    \delta^*_\epsilon \ge 4\pare{\frac{1}{\epsilon}-1} + \frac 23 (n-1) \Delta 
    \quad \in  \Omega\pare{\frac{1}{\epsilon} + n\Delta}.
\end{align*}
\end{proposition}

\begin{prf}
Let $n, \Delta, \epsilon$
be as in the statement, and consider problem \eqref{pr tight} with $a := \frac 12$,
$\beta := \frac 23$, 
and $t := \ceil{\frac{2}{\epsilon}-1} -1 \ge 0$.
Our assumptions imply $0 < a < (n-1)\beta\Delta.$ 
Therefore, Lemma~\ref{optial QP} implies that $x^d = -u$ is the unique optimal solution to \eqref{pr tight} and that $x^c = v$ is the unique optimal solution to the corresponding \eqref{pr QP}.

To prove the proposition, it suffices to show that the only $\epsilon$-approximate solution to \eqref{pr tight} is the optimal solution $-u$.
In fact, this implies $\delta^*_\epsilon = \norminf{x^c-x^d} = \norminf{v+u}$, and the latter norm can be bounded as follows:
\begin{align*}
    \norminf{v+u}
    & = 
    2 t + (n-1) \beta \Delta 
    =
    2\pare{\ceil{\frac{2}{\epsilon}-1}-1} + (n-1) \beta \Delta \\
    & \ge 
    2\pare{\frac{2}{\epsilon}-2} + (n-1) \beta \Delta
    =
    4\pare{\frac{1}{\epsilon}-1} + (n-1) \beta \Delta.
\end{align*}

Therefore, in the remainder of the proof we show that the only $\epsilon$-approximate solution to \eqref{pr tight} is the optimal solution $-u$.
If $\epsilon =1$, this is easy to see.
In fact, our definition of $t$ implies $t=0.$
Therefore, the unique feasible point for \eqref{pr tight} is the origin.
Therefore, in the remainder of the proof we assume $\epsilon \in (0,1)$.

Let $u':=(-t+1,0,\dots,0)$.
Note that $u'$ is in $\bar P$ since $\epsilon < 1$ implies $t \ge 1.$
We have
$f(u')=-\pare{t-\frac{1}{2}}^2=f(u).$
Furthermore, it is simple to see that any feasible vector for \eqref{pr tight} different from $\pm u,u'$ has cost strictly larger than $f(u)$.

It is simple to check that $f^d_{\max} = - \frac 14$, since the maximum is achieved at the origin.
The vector $u$ is an $\epsilon$-approximate solution to \eqref{pr tight} if and only if
\begin{align*}
\frac{f(u) - f(x^d)}{f^d_{\max} - f(x^d)} 
=
\frac{f(u)-f(-u)}{f^d_{\max}-f(-u)} 
= 
\frac{-\pare{t-\frac 12}^2 +  \pare{t+\frac 12}^2}{- \frac 1 4 +  \pare{t+\frac 12}^2} 
= 
\frac{2}{t+1} 
\le 
\epsilon.
\end{align*}

Note that our definition of $t$ implies $t < \frac 2\epsilon -1$, thus $\frac{2}{t+1} > \epsilon$.
This shows that the vector $u$ is not an $\epsilon$-approximate solution to \eqref{pr tight}. 
Since $-u$ is the only vector in $\bar P \cap \Z^n$ with cost strictly smaller than $f(u)$, the only $\epsilon$-approximate solution to \eqref{pr tight} is the optimal solution $-u$.
\end{prf}



Similarly to Proposition~\ref{prop IQP ub first}, also Proposition~\ref{prop IQP ub} implies that
if $k=n$, then
$\delta^*_\epsilon$ can grow at least linearly with respect to both $n$ and $\Delta$.
However, the bound in Proposition~\ref{prop IQP ub} is a function of $\epsilon$ as well.
It implies that $\delta^*_\epsilon$ can grow at least linearly with respect to $\frac 1 \epsilon.$

\subsubsection{Lower bounds on $\delta^\star_\epsilon$}
\label{sec tight b}

In this section we study the tightness of Theorem~\ref{th main}\ref{th main 2} by providing a lower bound on $\delta^\star_\epsilon$ that is a function of $n$, $\Delta$, and $\frac 1 \epsilon$.

\begin{proposition}
\label{prop new}
For every $n,\Delta \in \Z$ with $n\ge 2$, $\Delta \ge 2$ and $\epsilon \in (0,\frac 12)$, there exists an instance of \eqref{pr IQP} 
with subdeterminant $\Delta$ and $k=1$
for which
\begin{align*}
\delta^\star_\epsilon
\ge \frac{(n-1)\Delta-1}{\epsilon}-2 \quad \in \Omega \pare{\frac{n\Delta}{\epsilon}}.
\end{align*}
\end{proposition}

\begin{prf}
	Let $n,\Delta,\epsilon$ be as in the statement, and consider the \eqref{pr IQP}
	\begin{align}
	\label{pr tight in 8}
	\begin{split}
	\min \ & f(x)=
	-x_1^2 \\
	\stt \ & x \in \tilde P \cap \Z^n,
	\end{split}
	\end{align}
	where $\tilde P:= \bar P \cap \{x\in \R^n\mid x_1-\Delta\sum_{i=2}^nx_i\le \beta-1  +t \}.$ 
	Here, the polytope $\bar P$ is given in Definition~\ref{def poly} with $\beta := \frac 12$ and $t:=\floor{\frac{(n-1)\beta\Delta+\beta-1}{\epsilon}} \ge 1.$ 
Notice that the constraint matrix defining $\tilde P$ coincides with the one defining $\bar P$.
Therefore, Observation~\ref{obs delta} implies that problem \eqref{pr tight in 8} has subdeterminant $\Delta$.
	
Since $\tilde P \subseteq \bar P$, we have $\tilde P \cap \Z^n \subseteq \bar P \cap \Z^n$.
It can be easily checked that $-u \in \tilde P$, while $u \notin \tilde P$, therefore 
$x^d = -u$ is the unique optimal solution to \eqref{pr tight in 8}.

Let $w := ((n-1)\beta\Delta+\beta-1+t,\beta,\dots,\beta)$.
Observe that $w$ is a vector in $\tilde P$ with the largest $x_1$.
In fact, we know that every $x \in \bar P$ satisfies $x_i \le \beta$, for $i=2,\dots,n,$ thus for every $x \in \tilde P$ we have 
$x_1 \le \Delta\sum_{i=2}^{n}x_i + \beta-1+ t \le (n-1)\beta\Delta+\beta-1 + t$.
	Furthermore, the vector $-u$
	is a vector in $\tilde P$ with the smallest $x_1.$
In fact, we know that $-u$ is the vertex of $\bar P$ with the smallest $x_1$, and $-u \in \tilde P$.
	It can be checked that 
	\begin{align*}
		f(w) & = -((n-1)\beta\Delta+\beta-1+t)^2 < -t^2=f(-u),
	\end{align*}
	because $(n-1)\beta\Delta+\beta-1 > 2\beta-1=0. $
	In particular, we conclude that $x^c = w$ is an optimal solution to the \eqref{pr QP} corresponding to \eqref{pr tight in 8}.
	
	Next, we show that $-u$ is not an $\epsilon$-approximate solution to \eqref{pr QP}.
	Notice that $f^c_{\max}=0$ as it is achieved at the origin. 
	We have
	\begin{align*}
		\frac{f(-u)-f(x^c)}{f^c_{\max}-f(x^c)}
		&= 1-\frac{t^2}{((n-1)\beta\Delta+\beta-1+t)^2} \\
		& =1-\frac{1}{1+\frac{((n-1)\beta\Delta+\beta-1)^2}{t^2}+2\frac{(n-1)\beta\Delta+\beta-1}{t}} \\
		&\ge 1-\frac{1}{1+\epsilon^2+2\epsilon}=\epsilon\frac{2+\epsilon}{(1+\epsilon)^2}>\epsilon,
	\end{align*}
	where the first inequality holds because $t \le \frac{(n-1)\beta\Delta+\beta-1}{\epsilon},$ while the second inequality is correct because $\frac{2+\epsilon}{(1+\epsilon)^2} > 1,$ when $\epsilon \in (0,\frac 12).$
	Thus, $-u$ is not an $\epsilon$-approximate solution to \eqref{pr QP}.
	
Since $-u$ is a vector in $\tilde P$ with the smallest $x_1$, and due to the form of the objective function, every $\epsilon$-approximate solution $x^\star$ to \eqref{pr QP} must satisfy $x_1^\star > u_1 = t,$ which implies
	\begin{align*}
		\norminf{x^\star - x^d} > 2t \ge 2\frac{(n-1)\beta\Delta+\beta-1}{\epsilon}-2 = \frac{(n-1)\Delta-1}{\epsilon}-2.
	\end{align*}
	This implies $\delta^\star_\epsilon \ge \frac{(n-1)\Delta-1}{\epsilon}-2.$
\end{prf}

In particular, Proposition~\ref{prop new} shows that, in Theorem~\ref{th main}\ref{th main 2}, the dependence of $\delta^\star_\epsilon$ on $n$ is tight.
Furthermore, it shows that $\delta^\star_\epsilon$ can grow at least linearly with respect to both $\Delta$ and $\frac 1 \epsilon$.
Therefore, in Theorem~\ref{th main}\ref{th main 2}, the dependence on $\epsilon$ is tight for $k=1$.

	

\ifthenelse {\boolean{MPA}}
{
\bibliographystyle{spmpsci}
}
{
\bibliographystyle{plain}
}


\end{document}